\documentclass[]{article}

\usepackage{array}
\usepackage{amssymb}
\usepackage{graphicx}
\usepackage{epstopdf}
\usepackage{amsmath}
\usepackage{url}
\usepackage{bm}
\usepackage{graphicx}
\usepackage{caption}
\usepackage{multirow}
\usepackage{xcolor}
\usepackage{float}
\floatstyle{plaintop}
\restylefloat{table}
\begin{document}

\title{Efficient global optimization of constrained mixed variable problems}

\author{Julien Pelamatti, Lo\"ic Brevault, Mathieu Balesdent, \\ El-Ghazali Talbi, Yannick Guerin}
\date{}
\maketitle

\begin{abstract}
Due to the increasing demand for high performance and cost reduction within the framework of complex system design, numerical optimization of computationally costly problems is an increasingly popular topic in most engineering fields. In this paper,  several variants of the Efficient Global Optimization algorithm for costly constrained problems depending simultaneously on continuous decision variables as well as on quantitative and/or qualitative discrete design parameters are proposed. 
The adaptation that is considered is based on a redefinition of the Gaussian Process kernel as a product between the standard continuous kernel and a second kernel representing the covariance between the discrete variable values. Several parameterizations of this discrete kernel, with their respective strengths and weaknesses, are discussed in this paper. The novel algorithms are tested on a number of analytical test-cases and an aerospace related design problem, and it is shown that they require fewer function evaluations in order to converge towards the neighborhoods of the problem optima when compared to more commonly used optimization algorithms.
\end{abstract}

\section{Introduction}
Within the framework of complex system design, it is often necessary to solve optimization problems involving objective and constraint functions which depend simultaneously on continuous and discrete decision variables. These discrete variables can be either integer values (\emph{i.e.,} 1,2,3) or categorical design variables (\emph{e.g., materials, colors}). A common example of such an optimization problem is the preliminary design of aerospace vehicles, for which the performance criteria and constraints depend on continuous variables (\emph{e.g.,} sizing parameters) as well as on discrete design parameters (\emph{e.g.,} type of propulsion, type of materials). Furthermore, the estimation of the performance and constraint functions of such systems can be computationally costly, as for numerical simulations involving Finite Element Models (FEM) or coupled multidisciplinary analyses. Due to the complexity and computational cost of this type of problems, optimization algorithms commonly used in the presence of discrete variables such as  mixed variable Genetic Algorithm \cite{Stelmack1998} and Mesh Adaptive Discrete Search (MADS) \cite{Abramson2009}  result inadequate, as they require large number of function evaluations in order to conv
erge. It is therefore necessary to rely on optimization algorithms that can provide convergence towards the problem optimum with as few iterations as possible. An increasingly popular solution for computationally expensive problem is the Surrogate Model Based Design Optimization (SMBDO) \cite{Queipo2005}. SMBDO involves surrogate models of the numerical functions characterizing the problem, created by using a data set of limited size. These surrogate models are usually considerably cheaper to evaluate when compared to the exact problem functions. However, they also tend to introduce a modeling error which must be taken into account and dealt with. A simultaneous search for the problem optimum and refinement of the surrogate models is performed by first determining the location in the search space at which the problem optimum is most likely to be found according to a given infill criterion and by then calculating the actual system performance and constraint values at said location. Subsequently, the newly computed data sample is added to the pre-existing data set and the surrogate models are updated. This routine is repeated until a given convergence criterion is reached. In the literature, a number of different SMBDO algorithms have been proposed, each one relying on different Designs of Experiment (DoE), different surrogate modeling techniques and different infill criteria. A global overview of the most popular SMBDO techniques and their various aspects is provided by Queipo \emph{et. al.} \cite{Queipo2005}. More detailed reviews of surrogate modeling techniques for SMBDO purposes can be found in \cite{Haftka1998}, \cite{Simpson2001} and \cite{Wang2007a}. Finally, a review of the most popular data sample infill criteria is provided by Sasena \cite{Sasena2002}. 

A few surrogate modeling techniques for functions depending on both continuous and discrete variables exist in the literature \cite{Qian2008},\cite{Swiler2012},\cite{Zhang2015},  however, the majority of the proposed SMBDO techniques are defined to work within the continuous design space and only a small number has been extended to mixed continuous/discrete optimization problems. For instance, a number of variants of Radial Basis Functions (RBF) based SMDBO techniques for constrained mixed continuous/integer problems have been proposed \cite{Beauthier2014},\cite{Holmstrom2008},\cite{Muller2013},\cite{Rashid2008}. A surrogate model assisted GA is discussed in \cite{Bajer2013} where both Generalized Linear Models (GLM) and RBF are considered. Two main drawbacks can be identified within the existing mixed variable SMDBO techniques mentioned above. First of all, most of them are developed in order to deal with integer variables and can not handle optimization problems depending on generic discrete design variables (\emph{e.g.,} type of materials). Furthermore, the handling of the constraints in the previously mentioned optimization methods relies on direct penalization of the objective function values for solutions that are not feasible. Although popular, this approach results inadequate when confronted with expensive computations and usually requires large numbers of function evaluations.

In this paper, a mixed continuous/discrete adaptation of the Efficient Global Optimization (EGO), first proposed by Jones \cite{Jones1998}, is presented. EGO relies on Gaussian Process (GP) based surrogate models \cite{Rasmussen2006} and performs their refinement and optimization by using the Expected Improvement (EI) infill criterion. A number of possible extensions of EGO for the optimization of constrained problems is presented by Schonlau \cite{schonlau1998global}, notably the use of the Probability of Feasability (PoF), which is considered in this paper. In this work, an adaptation of the algorithm for the optimization of problems comprising both continuous and discrete decision variables is proposed and discussed. Following this introduction, in the second section the optimization problem is formally stated and the involved variables and parameters are defined. Subsequently, in the third section, the working principle of GP is briefly described and the considered kernels for the surrogate modeling of mixed variable functions are presented and compared. In Section 4, the proposed mixed variable adaptation of the EGO infill criterion and its optimization process are described and discussed. Afterwards, in the fifth section, these novel algorithms are tested on a number of analytical and aerospace design related benchmarks and their optimization performance is compared to commonly used algorithms. Finally, Section 6 presents the relevant conclusions which can be drawn from the presented results. Possible perspectives and improvements related to the proposed algorithms are also discussed.

\section{Problem statement}
\label{ProbStat}
Mixed-variable optimization problems involve continuous and discrete variables. Continuous variables (\emph{e.g.,} structure sizing parameters, combustion chamber pressure) refer to real numbers defined within a given interval. Discrete variables, instead, are here the non-relaxable variables defined within a finite set of choices. They often represent technological choices such as types of material and structural alternatives and can be either ordinal or nominal \cite{Agresti1996}. The first category includes variables for which a notion of order can be determined (\emph{e.g.,} integer variables, 'small/medium/large' types of decision variables), while the latter category includes unordered design parameters for which no concept of metrics can be defined (\emph{e.g.,} colors, materials).  Although discrete variables often lack a conceptual numerical representation, it is common practice to assign an integer value to every considered alternative in order to be able to include the related choices in the numerical optimization. Please note that the value assigned to a given discrete design choice is user-defined and has no influence on the presented optimization algorithm performance. For the sake of generality, in this paper no distinction between nominal and ordered discrete variables is made as the presented algorithms process both in the same way.
The considered mixed-variable optimization problem can be formulated as follows:
\begin{eqnarray}
\label{PBDefFirst}
\min &\;& f(\mathbf{x},\mathbf{z}) \qquad \qquad f: \mathcal{F}_x   \times \mathcal{F}_z  \rightarrow \mathcal{F}_y  \subseteq \mathbb{R}\\ 
\text{w.r.t.} &\;& \mathbf{x} \in \mathcal{F}_x \subseteq \mathbb{R}^{q} \nonumber \\
&& \mathbf{z}  \in \mathcal{F}_z  \nonumber \\
\text{s.t.} &\;& g_i  (\mathbf{x},\mathbf{z}) \leq  0 \qquad g_i: \mathcal{F}_x   \times \mathcal{F}_z  \rightarrow \mathcal{F}_{g_i}  \subseteq \mathbb{R} \label{PBDefLast}
 \\ 
&& \mbox{for } i = 1,...,n_g \nonumber
\end{eqnarray}

\noindent where $f(\cdot)$ is the  single-objective function characterizing the problem and defined in the codomain $\mathcal{F}_y$,  $g_i(\cdot)$ is one of the $n_g$  constraints the problem is subject to defined in the codomain $ \mathcal{F}_{g_i}$,  $\mathbf{x} = (x_1,...,x_q)$  is the vector containing the continuous decision variables and  $\mathbf{z} = (z_1,...,z_r)$  is the vector  containing the discrete decision variables. $q$ and $r$ being the size of the continuous and discrete dimensions of the functions $f(\cdot)$ and $g(\cdot)$. For the sake of simplicity, the input vector containing both continuous and discrete variables is represented with the following notation: $\mathbf{w} = \{ \mathbf{x}, \mathbf{z} \}$. Each discrete variable $z_j$ is characterized by $b_j$ possible values, also known as \emph{levels}, which therefore results in a total number of categorical combinations, or categories, of $m = \prod \limits_{k = 1}^{k = r} b_k$.

As mentioned in the introduction, the approach for the optimization of this kind of computationally intensive mixed variable problems that is considered and presented in this article relies on a mixed variable adaptation of EGO, which is discussed in the following sections of the paper.

\section{Gaussian process-based surrogate modeling of mixed variable functions}
\label{GaussMixedVar}
The core concept of surrogate modeling is to predict the response value $\hat{y}(\mathbf{w}^*)$ of a black-box function $f(\cdot)$ for a generic unmapped input $\mathbf{w}^* =  \{ \mathbf{x^*}, \mathbf{z^*} \}$ through an inductive procedure. The surrogate model is created by processing a so-called training data set $\mathcal{D}$, or Design of Experiments (DoE), of $n$ samples $\{\mathbf{w}^i,y^i\}$ with $i\in\{1,...,n\}$. $\mathcal{D}$ can be defined as follows:

\begin{equation*}
\mathcal{D} = \left\{ 
\mathcal{X} = \{ \mathbf{x}^1,..., \mathbf{x}^n\}, \quad 
\mathcal{Z} = \{ \mathbf{z}^1,..., \mathbf{z}^n\}, \quad 
\mathcal{Y} = \{ y^1,...,y^n\} 
 \right\}
\end{equation*} 

\noindent where $\mathcal{X}$ and $\mathcal{Z}$ are the subsets of $\mathcal{D}$ containing the continuous and discrete variables, respectively, while $\mathcal{Y}$ is the subset containing the corresponding responses. The surrogate modeling technique that is considered in this article is based on GP. A generic GP characterizes the probability distribution of the possible regression functions \cite{Rasmussen2006} and is usually defined by a mean $\mu(\mathbf{w})$ and a parameterized covariance function, or kernel, $k(\mathbf{w}^i,\mathbf{w}^j)$. The kernel is characterized by a number of hyperparameters, the optimal values of which can be determined through an internal optimization process known as training. The vast majority of covariance functions discussed in the literature are defined within the continuous domain and are invariant by translation (\emph{i.e.}, stationary), and the resulting surrogate model is often referred to as \emph{Kriging} \cite{Oliver1990},\cite{Sacks1989}. The prediction $\hat{y}$ of the function value at an unmapped location $\mathbf{w}^*$ is computed as the mean of the mean GP \textbf{$Y(\mathbf{w}^*)$} at said location :
\begin{equation}
\hat{y}(\mathbf{w}^*) =\overline{Y}(\mathbf{w}^*) = \mu +\bm{\psi}^T_{\mathcal{W}}(\mathbf{w}^*) \mathbf{K}^{-1}(\mathbf{y}-\mathbf{1}\mu) 
\label{Prediction}
\end{equation}
with:
\begin{equation}
\mu = \frac{\mathbf{1}^T\mathbf{K}^{-1}\mathbf{y}}{\mathbf{1}^T\mathbf{K}^{-1}\mathbf{1}}
\end{equation}
where $\mu$ is the GP mean that, in case no information regarding the modeled function
is known, is considered to be constant  \cite{Simpson2001} and $ \mathbf{K}$  is the $n\times n$ covariance matrix containing the covariance values between every sample of the data set:
\begin{equation}
\mathbf{K}_{i,j} = k(\mathbf{w}^i,\mathbf{w}^j)
\end{equation}
$\mathbf{y}$ is a   $n\times 1$ vector containing the responses corresponding to the $n$ data samples,  $\mathbf{1}$ is a $n\times 1$ vector of ones and finally $ \bm{\psi}_{\mathcal{W}} $ is an $n\times 1$ vector containing the covariance values between each sample of the training data set and the point at which the function is predicted:
\begin{equation}
 \bm{\psi}_{\mathcal{W}_i}(\mathbf{w}^*) =   k(\mathbf{w}^*,\mathbf{w}^i) \qquad \mbox{ for } i=1,...,n
 \label{CovVector}
\end{equation}
Furthermore, the estimate of the modeling uncertainty associated with a given prediction $\hat{y}(\mathbf{w}^*)$ can be computed under the form of a mean squared error as:
\begin{equation}
\hat{s}^2(\mathbf{w}^*) = k(\mathbf{w}^*,\mathbf{w}^*)  - \bm{\psi}^T_{\mathcal{W}}(\mathbf{w}^*) \mathbf{K
}^{-1} \bm{\psi}_{\mathcal{W}}(\mathbf{w}^*)   
\label{PredictionVar}
\end{equation}

The covariance function used to characterize a GP can be defined or adapted by the user. In order to ensure that the chosen kernel is valid, it is necessary for the covariance function to be symmetric and positive semi-definite \cite{Santner2003}. It is possible to show that the product between valid kernels also results in a valid one \cite{Shawe-Taylor2004}. This can be derived from the fact that according to the Schur product theorem, the Hadamard product between two positive semi-definite matrices results in a positive semi-definite matrix. In practice, the property that is used in this paper is that valid mixed variable kernels can be obtained by combining kernels defined in the continuous $q$-dimensional space and kernels defined in the discrete $r$-dimensional space, as is proposed by Roustant \emph{et al.} \cite{roustant:hal-01702607}:
\begin{equation}
k(\mathbf{w}^i,\mathbf{w}^j) = k_{c}(\mathbf{x}^i,\mathbf{x}^j)*k_{d}(\mathbf{z}^i,\mathbf{z}^j)
\end{equation}
where $ k_{c}$ and  $ k_{d}$ represent the continuous and discrete kernels, respectively, while $*$ is a generic operator allowing to combine kernels. Examples of such an operator found in the literature are product, sum and ANOVA. In this article, only the kernel-wise product is considered. It is important to note that no assumption on the type of input data is required in order to define a valid kernel function, as long as it can be associated to a probability measure. By extension, kernels computed on continuous and discrete variables can be combined without any loss of generality or applicability in order to obtain the covariance matrix of a multivariate Gaussian variable. In the most generic case, the resulting kernel is considered to be heteroscedastic, \emph{i.e.,} characterized by a different GP variance for each discrete category of the problem. 

A large number of valid continuous kernels exist in the literature, however, depending on the characteristics of the function that is being modeled, some choices result more suitable than others. In this paper, the \emph{p-exponential} covariance function \cite{Santner2003} is considered:
\begin{equation}
k_{c}(\mathbf{x}^i,\mathbf{x}^j) = \sigma_c^2 \exp\left(- \sum _{k = 1}^{k = q} \theta_k |x_k^i-x_k^j|^{p_k} \right)
\end{equation}

\noindent where $x_k^i$ is the value of the continuous variable corresponding to the $k^{th}$ dimension of the data sample $\mathbf{x}^i$,  $\theta_k$ and $p_k$ are hyperparameters characterizing the surrogate model and $\sigma_c^2$ is the variance associated to the continuous kernel.

\subsection{Discrete kernel complete parameterization }
\label{CompleteParam}
In the most generic case, the discrete kernel $k_d$ can be represented under the form of a $m \times m$ positive semi-definite matrix $\mathbf{T}$. Each element $\mathbf{T}_{k,s}$ of this matrix contains the covariance between two generic discrete categories $k$ and $s$ of the modeled function. A proper parameterization of the matrix $\mathbf{T}$ is necessary in order to ensure the validity of the kernel. In this paper, the Cholesky decomposition of the covariance matrix \cite{Pinheiro1996} is considered:
\begin{equation}
\mathbf{T} = \mathbf{L}\mathbf{L}^T
\end{equation}
in which $\mathbf{L}$ is a lower triangular matrix. The matrix $\mathbf{L}$ is created with the help of the so-called hypersphere decomposition \cite{Rebonato2011}, first applied within the scope of mixed variable GP by Zhou \emph{et al.} \cite{Zhou2011}. The underlying idea is that the elements of the $k^{th}$ row of the matrix represent the coordinates of a point on the surface of a $k$-dimensional hypersphere. The triangular matrix elements $l_{k,s}$ are defined as follows:
\begin{equation}
	\left\{
                \begin{array}{ll}
		l_{1,1} = \alpha_{ k,0} \\
		l_{k,1} = \alpha_{ k,0} \cos(\alpha_{ k,1}) \\
		l_{k,s} = \alpha_{ k,0} \sin(\alpha_{ k,1}) ... \sin(\alpha_{ k,s-1}) \cos(\alpha_{ k,s})
		\qquad \mbox{for } s = 2,...,k-1 \\
		l_{k,k} = \alpha_{ k,0} \sin(\alpha_{ k,1}) ... \sin(\alpha_{ k,k-2}) \sin(\alpha_{ k,k-1})		
                \end{array}
              \right.
              \label{L}
\end{equation}

\noindent where $\alpha_{k,0} > 0 $ and $\alpha_{k,s} \in (0,\pi)$ (for $s \neq 0$) are the hyperparameters characterizing the covariance between the various discrete categories of the modeled function. Although there is a total of $m^2$ combinations of discrete categories, due to the inherent symmetrical  nature of the covariance, $(m+1)m/2$ parameters $\alpha_{k,s}$ are required in order to define the matrix  $\mathbf{L}$. For complex problems characterized by a large number of discrete categories, the required number of hyperparameters becomes considerably large, thus making it very problematic to determine their optimal value \cite{Zhou2011}. Furthermore, in order to perform the optimization process of said hyperparameters, it is necessary for the data set to contain samples belonging to every discrete category of the considered problem. However, when dealing with actual engineering design problems, it may occur that some combinations of discrete design variables are not physically feasible or can not be modeled, in which case the previously complete parameterization of the discrete kernel can not be applied. An example of this issue can be found in the modeling of a rocket engine propulsive performance, in which not all the combinations of reductant and oxidant result in a feasible combustion process \cite{Pelamatti2018a}. For these reasons, the complete parameterization of the discrete kernel is not considered in this paper. In order to avoid these issues, alternative parameterizations which do not require all the discrete categories to be present within the data set and which are characterized by smaller numbers of hyperparameters are discussed in the following paragraphs.

\subsection{Alternative discrete kernel parameterizations}
\label{AlternParam}
In order to reduce the number of hyperparameters required to characterize the covariance matrix, a choice that can be made is to define the discrete kernel as a combination of 1-dimensional discrete kernels by relying on the operations listed in the previous paragraphs. For the sake of clarity, each 1-dimensional kernel can be re-written under the form of a positive semi-definite matrix. In case a kernel-wise product is considered, the resulting discrete kernel can be defined as \cite{roustant:hal-01702607}:
\begin{equation}
k_{d}(\mathbf{z}^i,\mathbf{z}^j) = \prod_{s = 1}^{r} [\mathbf{T}_{s}]_{z_s^i,z_s^j}
\end{equation}
where each matrix $\mathbf{T}_{s}$ contains the values of the covariance between the various levels of the generic discrete variable $s$. In the following paragraphs, possible parameterizations of $\mathbf{T}_{s}$ are described and discussed.

\subsubsection{Heteroscedastic dimension-wise hypersphere decomposition}
A first possible parametrization for  $\mathbf{T}_{s}$ that is considered, suggested by Zhou \cite{Zhou2011}, relies on the hypersphere decomposition described in the previous paragraph. In this case, the covariance is characterized dimension-wise rather than category-wise, the number of hyperparameters required to define the $r$ matrices $\mathbf{T}_{s}$ is therefore equal to $\sum_{k = 1}^{k = r} b_k(b_k+1)/2$. 
When compared to the complete hypersphere decomposition of the discrete kernel presented in Section \ref{CompleteParam}, the dimension-wise variant offers a better scaling with the discrete dimension of the problem that is being modeled in terms of number of hyperparameters, but as a trade-off provides a theoretically less accurate modeling of the correlation between the various discrete categories of the problem \cite{Pelamatti2018a}. Furthermore, due to the fact that the kernel is defined dimension-wise, it is not necessary for all the problem categories to be represented in the training data set in order to train the hyperparameters.

\subsubsection{Homoscedastic dimension-wise hypersphere decomposition}
If a further reduction of the number of hyperparameters is required, the assumption of homoscedasticity can be made. In other words, it can be assumed that all the categories of the problem are characterized by the same variance value. In this case, each matrix $\mathbf{T}_{s}$ has a constant diagonal value. By consequence, the discrete kernel can be rewritten as a product between the common GP variance and $r$ dimension-wise correlation matrices parameterized with the help of the same hypersphere decomposition discussed in the previous paragraphs. Due to the fact that correlation matrices are characterized by a unit diagonal, their hypersphere decomposition only requires  $\sum_{k = 1}^{k = r} b_k(b_k-1)/2$ hyperparameters (\emph{i.e.,} all the hyperparameters $\alpha_{i,0}$ are equal to 1).
Nevertheless, it must be kept in mind that the assumption of homoscedasticity can introduce a considerable modeling error when dealing with optimization problems that present categories characterized by considerably different behaviors, and might therefore not always be valid.
\subsubsection{Compound symmetry parameterization}
In case an ulterior reduction in the number of hyperparameters characterizing the discrete kernel is required, a parameterization that can be considered relies on the so-called Compound Symmetry (CS) \cite{Pinheiro2009}. In this case, each matrix $\mathbf{T}_s$ is characterized by a single value of covariance $c_s$ and a single value of variance $v_s$:
\begin{equation}
\label{CS_Matrix}
[\mathbf{T}_{s}]_{z_s^i,z_s^j} =  \begin{cases} 
v_s \quad \mbox{if} \quad z_s^i = z_s^j  \\ 
c_s  \quad \mbox{if} \quad z_s^i \neq z_s^j 
\end{cases}
\end{equation}
with $-(b_s + 1 )^{-1} v_s < c_s < v_s$, in order to ensure that $\mathbf{T}_{s}$ is positive semi-definite. A particular case of CS parameterization that is discussed in this article can be obtained by considering the covariance in the mixed continuous discrete search space to be spatially dependent as a function of the so-called Gower distance, as is proposed by Halstrup \cite{Halstrup2016}. In the Gower distance \cite{Gower1971}, the coordinates of the two samples that are being considered are compared dimension-wise. For the continuous dimensions, the distance is proportional to the Manhattan distance, while for the discrete dimensions the distance is a weighted binary value which depends on the similarity between the variable values. In practice, the Gower distance between two points can be expressed as follows:

\begin{equation}
d_{gow} (\mathbf{w}^i, \mathbf{w}^j) = \frac{\sum_{k = 1}^{k = q}  \frac{| x_{k}^i -  x_{k}^j|} {\Delta x_k}}{r+q} + \frac{\sum_{k = 1}^{k = r} S( z_{k}^i ,  z_{k}^j)  }{r+q}
\end{equation}

\noindent where $\Delta x_k$ is the range of the continuous variable in the $k$-dimension and $S$ is a score function defined as:
\begin{equation}
S( z_{k}^i ,  z_{k}^j) = \left\{
                \begin{array}{ll}
		0 \mbox{ if } z_{k}^i =  z_{k}^j \\ 
		1 \mbox{ if } z_{i}^i \neq  z_{k}^j
		
                \end{array}
              \right.
\end{equation}

\noindent Having re-defined the distance in order to take into account the presence of  both continuous and discrete variables, the \emph{p-exponential} covariance function may be used in order to define the mixed variable kernel:
\begin{equation}
k(\mathbf{w}^i,\mathbf{w}^j) = \sigma^2 \exp \left[-  \sum_{k = 1}^{k = q} \theta_k \left( \frac{  \frac{| x_{k}^i -  x_{k}^j|} {\Delta x_k}}{r+q} \right)^{p_k} - \sum_{k = 1}^{k = r} \theta_{k+q} \left( \frac{  S( z_{k}^i ,  z_{k}^j)}{r+q} \right)^{p_{k+q}} \right]
\end{equation}
This is equivalent to defining the mixed variable kernel as:
\begin{equation}
k(\mathbf{w}^i,\mathbf{w}^j) = k_c(\mathbf{x}^i,\mathbf{x}^j)  * \prod_{s = 1}^{r} [\mathbf{T}_{s}]_{z_s^i,z_s^j}
\end{equation}
where each matrix $\mathbf{T}_{s}$ is a CS covariance matrix defined as shown in Eq.(\ref{CS_Matrix}) with a $c_s/v_s$ ratio equal to:
\begin{equation}
c_s/v_s = \exp  \left[ -  \theta_{s} \left( \frac{  S( z_{s}^i ,  z_{s}^j)}{r+q} \right)^{p_{s}} \right] 
\end{equation} 
This CS Gower distance based GP is characterized by $ 2(q+r)$  hyperparameters and scales therefore better with the discrete dimension of the problem when compared to the previously described kernel parameterizations. Furthermore, the adaptation of a standard Gaussian process into a Gower distance based one is relatively simple. However, it must also be noticed that because each discrete variable is only parameterized by 2 hyperparameters $\theta$ and $p$, the surrogate model may present poor modeling performances in the presence of discrete variables with large numbers of discrete levels. For the same reason, the simultaneous presence of correlation and anti-correlation trends described by the same discrete variable might be more difficult to model when compared to the parameterizations described in the previous paragraphs \cite{Pelamatti2018a}.

\subsubsection{Discrete kernel parameterization examples}
As an illustrative example, the matrices characterizing the discrete kernel for a function characterized by 2 discrete variables with 3 levels each ($z_1 \in \{0,1,2\}$, $z_2 \in \{0,1,2\}$) resulting in 9 categories are provided below. In the heteroscedastic hypersphere decomposition case, the lower triangular Cholesky decomposition matrices used for the computation of $\mathbf{T}_s$ (with $s = 1,2$), depend on 6 hyperparameters $\alpha$ each, for a total of 12 hyperparameters. The generic matrix $\mathbf{L}_s$ has the following expression:
\begin{eqnarray}
\small
              \mathbf{L}_s =  \left[ \begin{array}{c |c |c }
                \alpha_{{1,0}_s} & 0 & 0\\ \hline
                \alpha_{{ 2,0}_s} \cos(\alpha_{{ 2,1}_s})  & \alpha_{{ 2,0}_s} \sin(\alpha_{{ 2,1}_s}) & 0 \\ \hline
                \alpha_{{ 3,0}_s} \cos(\alpha_{{ 3,1}_s})  & \alpha_{{ 3,0}_s} \sin(\alpha_{{ 3,1}_s}) \cos(\alpha_{{ 3,2}_s}) & \alpha_{{3,0}_s} \sin(\alpha_{{ 3,1}_s})\sin(\alpha_{{ 3,2}_s}) 
                \end{array}
                \right] 
               \label{2L4}
\end{eqnarray}
If the homoscedastic case is considered, the hyperparameters characterizing the category-specific variance are not necessary. By consequence, the matrix $\mathbf{L}_s$ only requires 3  hyperparameters $\alpha$ in order to be defined, thus resulting in a total of 6 hyperparameters, and acquires the following expression:
\begin{eqnarray}
\small
              \mathbf{L}_s =  \left[ \begin{array}{c |c |c }
                1 & 0 & 0\\ \hline
                 \cos(\alpha_{{ 2,1}_s})  & \sin(\alpha_{{ 2,1}_s}) & 0 \\ \hline
                 \cos(\alpha_{{ 3,1}_s})  &  \sin(\alpha_{{ 3,1}_s}) \cos(\alpha_{{ 3,2}_s}) &  \sin(\alpha_{{ 3,1}_s})\sin(\alpha_{{ 3,2}_s}) 
                \end{array}
                \right] 
\end{eqnarray}
Finally, if the compound symmetry decomposition of the discrete kernel is considered, no intermediate lower triangular matrix is necessary, and the matrix $\mathbf{T}_s$ can be computed directly as:
\begin{eqnarray}
\small
              \mathbf{T}_s =  \left[ \begin{array}{c |c |c }
                v_s & c_s & c_s\\ \hline
                 c_s  & v_s & c_s\\ \hline
                 c_s  &  c_s &  v_s 
                \end{array}
                \right] 
\end{eqnarray}
In this case, each matrix $\mathbf{T}_s$ is characterized by 2 hyperparameters, for a total of 4 hyperparameters. 
\section{Efficient global optimization of mixed variable functions}
As mentioned in the previous sections, the approach that is presented in this article relies on an adaptation of the EGO algorithm \cite{Jones1998} for the optimization of mixed variable constrained problems. The EGO routine can be divided in two main phases: the first phase consists in creating a GP based surrogate model of the objective and constraint functions from an initial data set. Subsequently, additional data samples are evaluated and added to the data set with the purpose of simultaneously refining the surrogate model and exploring the areas of the search space more likely to contain the optimization problem optimum. This refinement process is often referred to as infill. The location at which the newly added data samples are computed is determined according to a given infill criterion. In this article, the infill criterion to be considered is the Expected Improvement (EI) combined with the Probability of Feasibility (PoF). Please note that alternative solutions for the handling of constraints, such as the Expected Violation (EV) and the Constrained Expected Improvement (CEI) exist. The PoF is selected due to its robustness and constraint prediction accuracy \cite{Durantin2016}.
The contribution that is proposed in this article consists in performing the SMBDO of constrained mixed variable problems by relying on mixed continuous/discrete GP. In order to perform this sort of SMBDO, the data sample infill criterion must also be adapted to the presence of discrete variables. In the following sections, the used mixed variable infill criterion is described and commented. Subsequently, the approach and necessary considerations for the optimization of said criterion are briefly discussed.

\subsection{Mixed variable infill criterion}
The mixed variable Infill Criterion (IC) that is used in this article is a mixed variable adaptation of the one discussed by  Schonlau \cite{schonlau1998global}, defined as the product between EI and PoF. As is mentioned in Section \ref{GaussMixedVar}, the mixed variable GP kernel is defined in such a way that the resulting covariance matrix can be used to characterize a normal distribution. By extension, the derivation of the EI and the PoF expression remain valid when applied in the mixed variable search space rather than in a purely continuous one.
EI represents the expected value of the predicted improvement with respect to the data set:
\small
\begin{eqnarray}
\hspace{-0.7cm}
\mathbb{E}[I(\mathbf{w^*})] & = & \mathbb{E} \left[\max \left(y_{min} - Y(\mathbf{w}^*),0 \right) \right] \\
 & = & (y_{min} - \hat{y}(\mathbf{w^*}))\Phi \left(\frac{y_{min}-\hat{y}(\mathbf{w^*})}{\hat{s}(\mathbf{w^*})}\right) + \hat{s}(\mathbf{w^*})\phi \left(\frac{y_{min}-\hat{y}(\mathbf{w^*})}{\hat{s}(
\mathbf{w^*})}\right)
\label{EI}
\end{eqnarray}
\normalsize
\noindent where $y_{min}$ is the current minimum value present within the data set, while $\Phi(.)$ and $\phi(.)$  are the standard distribution and normal density functions, respectively. PoF, instead, represents the probability that all the constraints the problem is subject to are met at the unmapped location $\mathbf{w}^*$ of the search space. Given a constraint function $g_i(.)$, the probability for it to be satisfied at $\mathbf{w}^*$ can be computed as:

\begin{equation}
\mathbb{P}(g_i(\mathbf{w}^*) \leq 0) =  \Phi\left(\frac{0-\hat{g}_i(\mathbf{w}^*)}{\hat{s}_{g_i}(\mathbf{w}^*)}\right) 
\label{PF}
\end{equation}

\noindent where $\hat{s}_{g_i}$ refers to the estimated error in the prediction $\hat{g}_i$ of the constraint function. By extension, the PoF for problems subject to $n_g$ constraints can be computed as:

\begin{equation}
PoF ( \mathbf{w}^* )= \prod_{i = 1}^{n_g} \mathbb{P}(g_i(\mathbf{w}^*)\leq 0)
\label{PF2}
\end{equation}
The constrained optimization infill criterion $IC$ used in this article can finally be computed as:
\begin{equation}
IC ( \mathbf{w} ^*) = \mathbb{E}[I(\mathbf{w^*})] PoF ( \mathbf{w}^* )
\label{IC}
\end{equation}
The data sample to be added to the GP training data set is obtained by evaluating the value of the objective and constraint functions for the value of $\mathbf{w}^*$ that maximizes the IC:

\begin{equation}
\mathbf{w^*} = \mbox{argmax} (IC(\mathbf{w}))
\end{equation}
A secondary optimization process is therefore necessary. However, given that the computation time required to evaluate the IC is negligible when compared to the computation of the objective and constraint functions involved in complex system design problems, more common optimization algorithms can be used, as is discussed in the following part of the section. 
Once the value of $\mathbf{w}^*$ that maximizes IC has been determined, the exact objective and constraint functions of the optimization problem are computed at said location and the obtained data sample is added to the GP data set. Subsequently, the surrogate model must be trained anew in order to take into account the additional information provided by the latest data sample. This process is repeated until a user-defined stopping criterion is reached. For instance, in the original formulation of EGO the optimization is considered to have converged when the maximum EI value is smaller than 1\% of the best current function value. For the sake of simplicity and reproducibility of results, in this article an initial computational budget is defined and the optimization processes are stopped once the predefined number of data samples has been infilled.

\subsection{Infill criterion optimization}
As previously mentioned, in order to determine the location in the search space at which the infill data sample must be computed, the IC must be optimized. Said optimization can either be performed separately in each category of the problem by subsequently choosing the category yielding the largest value or it can be directly performed in the mixed continuous/discrete search space. In this article, the latter option is chosen. In this case, the IC is defined in a mixed variable search space (the same as the optimization problem) and by consequence most of the commonly used algorithms, such as gradient based ones, may not be used. Furthermore, due to the non-linear nature of the objective function and to the presence of several local optima, an evolutionary algorithm is used. More specifically, the IC is optimized by relying on a mixed continuous/discrete Genetic Algorithm (GA) similar to the one presented by Stelmack \cite{Stelmack1998} and coded with the help of the python based toolbox DEAP \cite{DEAP_JMLR2012}. The GA optimization routine is terminated once the maximum number of generations has been reached, or alternatively after the optimum value of the objective function has not improved over a predefined number of generations.

\section{Applications and Results}
In this section, the performance of the mixed variable adaptations of EGO based on the various mixed variable kernel definitions described in Section \ref{AlternParam}, \emph{i.e.,} heteroscedastic and homoscedastic dimension-wise hypersphere decomposition (He$\_$HS and Ho$\_$HS, respectively) as well as Gower distance based CS decomposition (CS) are compared on a number of analytical and aerospace design related optimization problems. In order to provide a measure of comparison with optimization methods commonly used, the results obtained with a penalized  mixed variable GA \cite{Stelmack1998} and an EGO based on independent GP for each category of the considered problem, or Category-Wise EGO (CW), are presented as well. In the first part of the section, results obtained on analytical test-cases are detailed. More specifically, mixed variable adaptations of the Branin and Goldstein functions characterized by varying continuous and discrete dimensions are considered. In the second part of the section, results obtained for the optimization of a launch vehicle booster performance under a number of feasibility constraints are discussed.

For all the results presented in this Section, the continuous part of the initial training data set is sampled with the help of a stochastic Latin Hypercube Sampling (LHS) \cite{McKay1979}. Subsequently, an even number of data samples is randomly assigned to each category characterizing the considered problem. The same initial data set is used for all the SMBDO methods in order to facilitate the comparison of results. Furthermore, in order to quantify and compensate the influence of the initial DoE random nature, each optimization problem is solved 10 times with different initial training data sets. The considered surrogate model hyperparameters are trained through the optimization of the likelihood function \cite{Rasmussen2006} by relying on the Covariance Matrix Adaptation Evolution Strategy (CMA-ES) \cite{Hansen2006}. In the homoscedastic cases, the optimal value of the GP variance is determined analytically as:
\begin{equation}
\sigma^2 = \frac{ (\mathbf{y} - \mathbf{1}\mu)^T\mathbf{R}^{-1} (\mathbf{y} - \mathbf{1}\mu)}{n}
\end{equation}
where $\mathbf{R}$ is an $n \times n$ matrix containing the values of the correlations between the data samples of the training set.

\subsection{Branin function}
\begin{figure}[h!]
\centering
\begin{minipage}{.5\textwidth}
  \centering
  \includegraphics[width=1.1\linewidth]{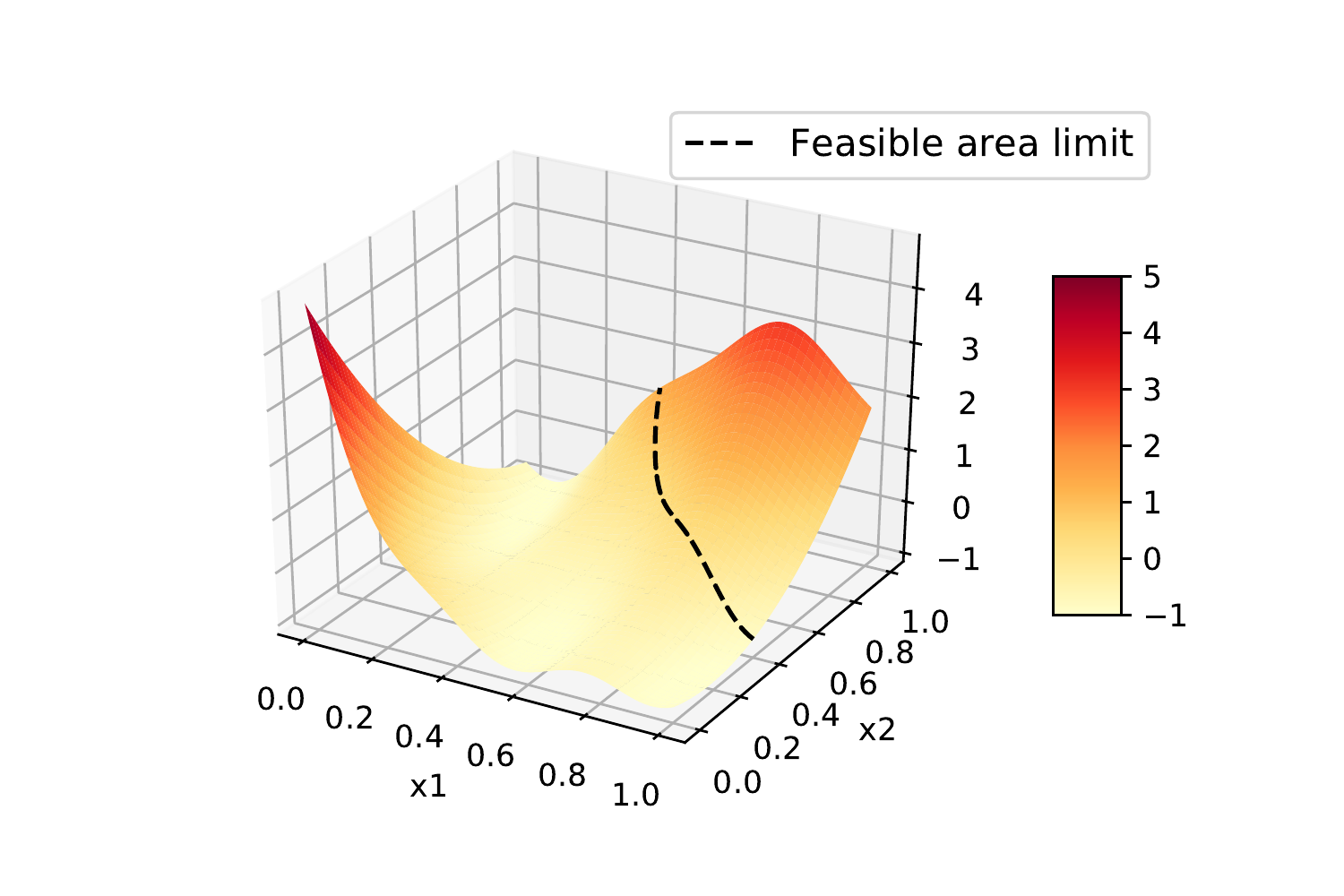}
  \captionsetup{labelformat=empty}
  \caption*{Category 1: $z_1 = 0, z_2 = 0$}
\end{minipage}%
\begin{minipage}{.5\textwidth}
  \centering
  \includegraphics[width=1.1\linewidth]{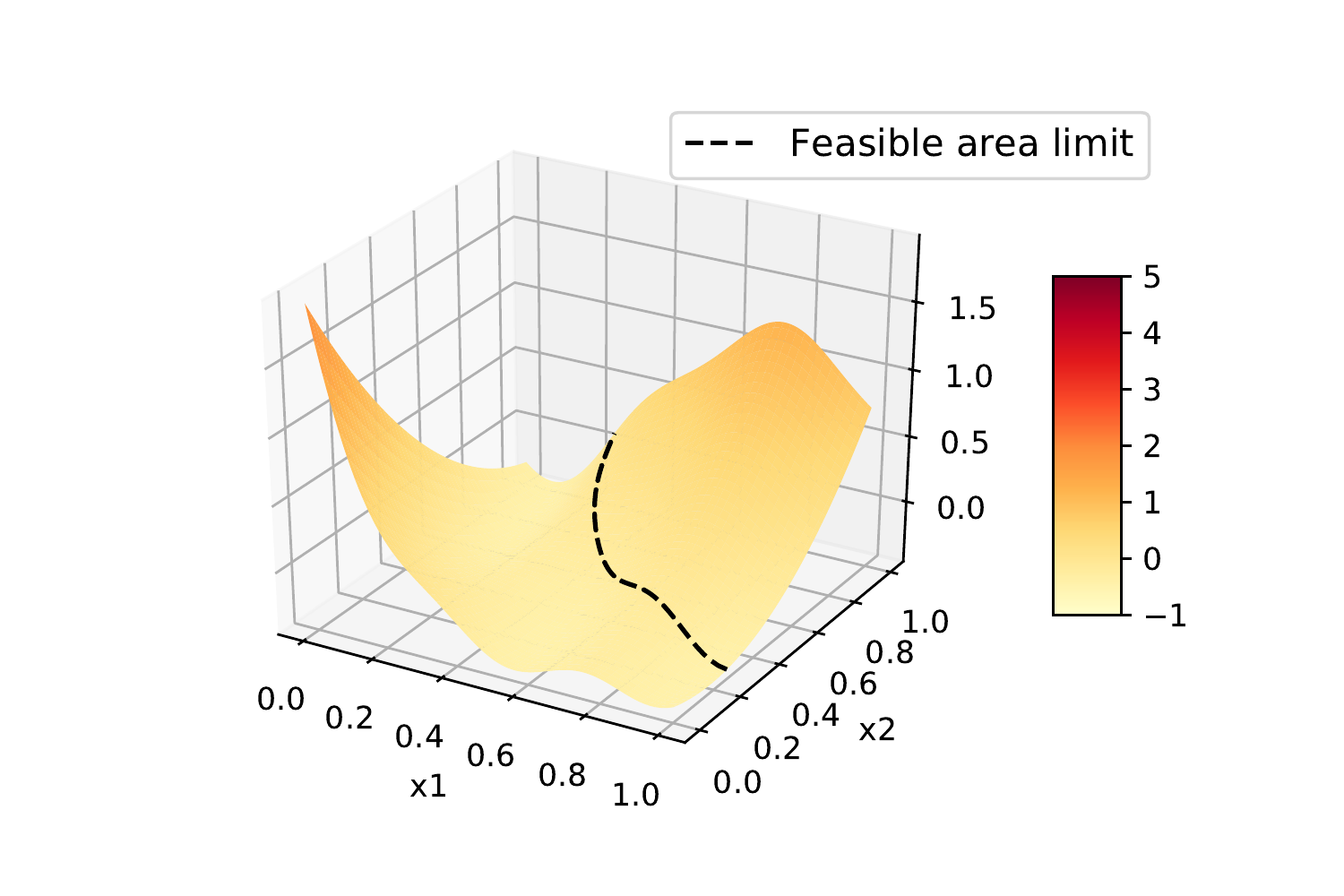}
  \captionsetup{labelformat=empty}
  \caption*{Category 2: $z_1 = 0, z_2 = 1$}
\end{minipage}
\begin{minipage}{.5\textwidth}
  \centering
  \includegraphics[width=1.1\linewidth]{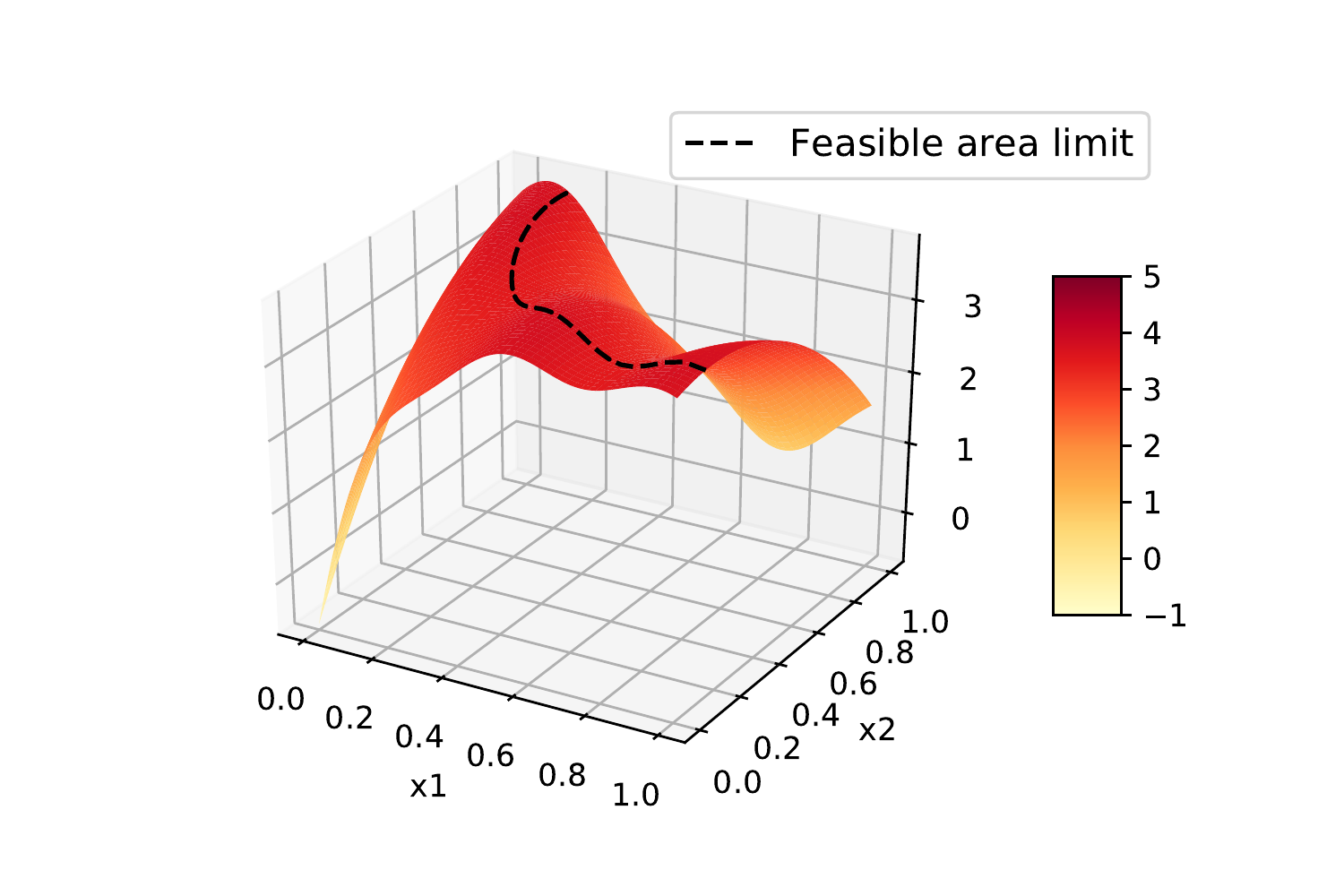}
  \captionsetup{labelformat=empty}
  \caption*{Category 3: $z_1 = 1, z_2 = 0$}
\end{minipage}%
\begin{minipage}{.5\textwidth}
  \centering
  \includegraphics[width=1.1\linewidth]{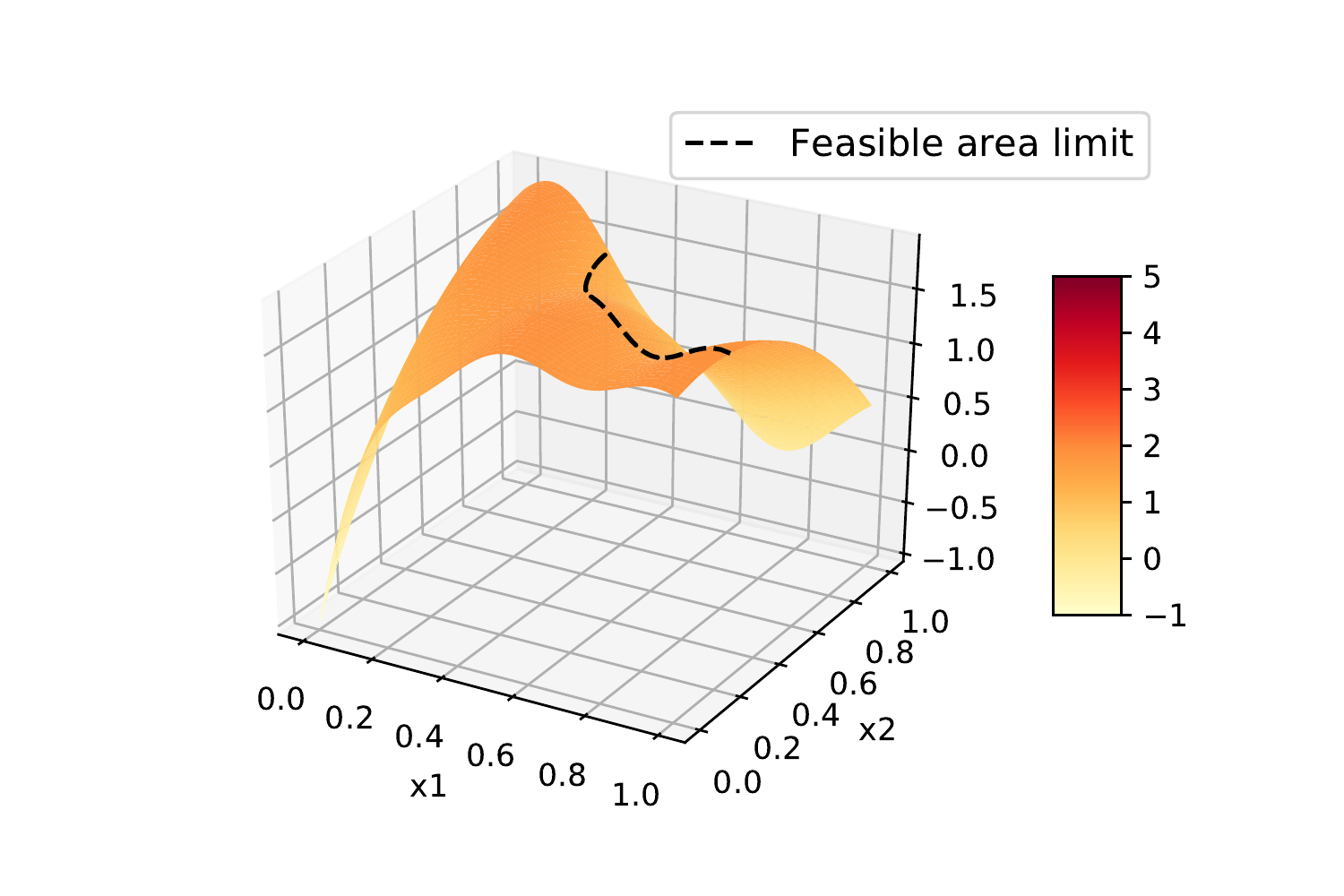}
  \captionsetup{labelformat=empty}
  \caption*{Category 4: $z_1 = 1, z_2 = 1$}
\end{minipage}
\caption{The 4 discrete categories of the Branin function}
\label{Branin}
\end{figure}

\begin{figure}[h!]
 \includegraphics[width=0.9\linewidth]{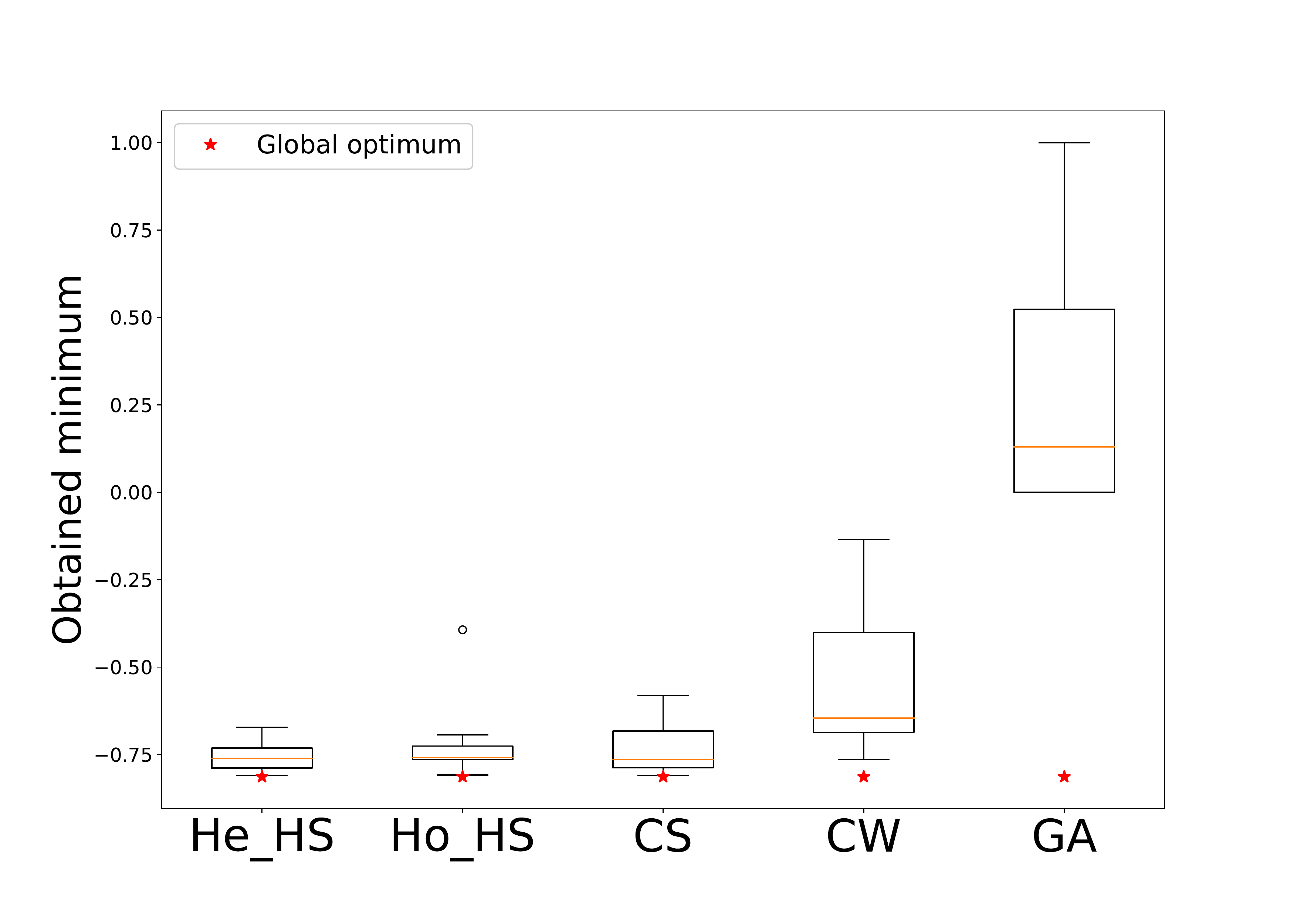}
  \caption{Optimization results for a constrained Branin function obtained over 10 repetitions. From left to right, the used optimization techniques rely on heteroscedastic dimension-wise decomposition (He$\_$HS),  homoscedastic dimension-wise decomposition (Ho$\_$HS), compound symmetry decomposition (CS), category-wise separate surrogate modeling (CW) and finally the last results are obtained with the help of a Genetic Algorithm (GA)}
  \label{BraninRes1}
  \end{figure}

\begin{figure}[h!]
 \includegraphics[width=0.9\linewidth]{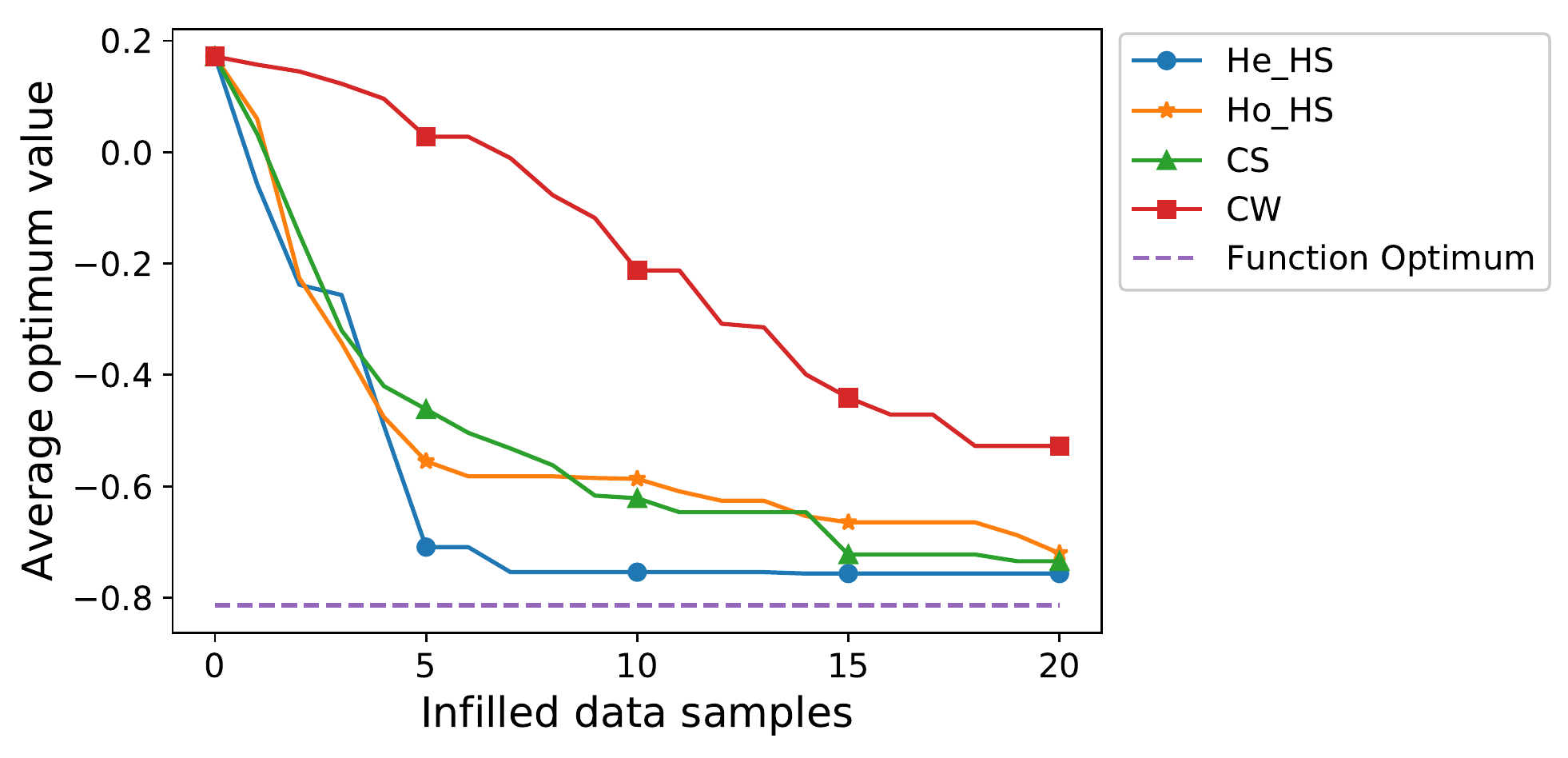}
  \caption{Evolution of the obtained optimum value for the optimization of the Branin test-case as a function of the number of infilled data samples averaged over 10 repetitions. The used optimization techniques rely on heteroscedastic dimension-wise decomposition (He$\_$HS),  homoscedastic dimension-wise decomposition (Ho$\_$HS), compound symmetry decomposition (CS) and category-wise separate surrogate modeling (CW)}
  \label{BraninRes2}
  \end{figure}

The first analytical benchmark to be considered is a modified version of the Branin function characterized by two continuous variables and two discrete variables, each one presenting 2 levels, thus resulting in a total of 4 discrete categories. This adaptation of the Branin function optimization problem is defined as follows:

\begin{eqnarray}
&& \min  f(x_1,x_2,z_1,z_2) \\ 
&& \text{w.r.t.} \quad x_1,x_2,z_1,z_2 \nonumber
\end{eqnarray}
\begin{equation}
\qquad \text{s.t.} \quad g(x_1,x_2,z_1,z_2) \geq 0
\end{equation}
with:
\begin{eqnarray*}
x_1 \in [0,1], \quad x_2 \in [0,1], \quad z_1 \in \{0,1\}, \quad z_2 \in \{0,1\} 
\end{eqnarray*}

\noindent where:
\begin{equation}
  f(x_1,x_2,z_1,z_2)=\begin{cases}
               h(x_1,x_2) \hfill   \qquad  \mbox{ if }  z_1 = 0  \mbox{ and } z_2 = 0\\
               0.4 h(x_1,x_2) \hfill  \qquad  \mbox{ if } z_1 = 0 \mbox{ and } z_2 = 1\\
               -0.75h(x_1,x_2) +3.0  \hfill  \qquad  \mbox{ if } z_1 = 1 \mbox{ and } z_2 = 0 \\
               -0.5h(x_1,x_2)  +1.4 \hfill  \qquad  \mbox{ if } z_1 = 1 \mbox{ and } z_2 = 1
            \end{cases}
\end{equation}

\begin{equation}
\hspace{-1.5cm}
\begin{aligned}
h(x_1,x_2) = & \left[ \left( (15x_2-\frac{5}{4\pi^2}(15x_1-5)^2+\frac{5}{\pi}(15x_1-5)-6)^2 + \right. \right.
\\ & \left. \left. 10\left(1-\frac{1}{8\pi}\right)\cos(15x_1-5)+10 \right) -54.8104 \right] \frac{1}{51.9496}
\end{aligned}
\end{equation}

\begin{equation}
  g(x_1,x_2,z_1,z_2)=\begin{cases}
               x_1x_2-0.4 \hfill   \qquad  \mbox{ if }  z_1 = 0  \mbox{ and } z_2 = 0\\
               1.5x_1x_2-0.4 \hfill  \qquad  \mbox{ if } z_1 = 0 \mbox{ and } z_2 = 1\\
              1.5x_1x_2-0.2 \hfill  \qquad  \mbox{ if } z_1 = 1 \mbox{ and } z_2 = 0 \\
               1.2x_1x_2-0.3 \hfill  \qquad  \mbox{ if } z_1 = 1 \mbox{ and } z_2 = 1
            \end{cases}
\end{equation}
The 4 categories of the considered Branin function are presented in Figure \ref{Branin}.
For the SMBDO techniques, an initial training data set of 20 samples is used and subsequently 20 additional data points are infilled during the optimization process. Instead, the GA is initialized with a population of 5 individuals which evolve for 8 generations, thus resulting in the same number of exact function evaluations as the SMBDO techniques. The results obtained for the optimization of the previously described constrained Branin function over 10 repetitions are presented in Figures \ref{BraninRes1}  and \ref{BraninRes2} as well as in Table \ref{TableBranin}. Please note that in Table \ref{TableBranin} the average constraint value at convergence is not provided for the GA optimizations due to the fact that it relies on a penalized objective function.
\begin{table}
\small
\begin{tabular}{|c|c|c|c|c|}
\hline
\multirow{ 2}{*}{Method} & Average & Average   & \# of optima in & \# of \\ 
 & optimum value & constraint value  &  the correct category & hyperparameters \\ \hline
 He$\_$HS & -0.689  $\pm$  16.51\% & 0.026 & 9 & 10 \\ \hline
Ho$\_$HS & -0.784  $\pm$ 4.78\% & 0.009 & 10 & 6 \\ \hline
CS & -0.799  $\pm$ 1.53\% & 0.008 & 10 & 8 \\ \hline
CW &  -0.596  $\pm$ 25.85\% & 0.043 & 7 & 4/category \\ \hline
GA & -0.158  $\pm$  56.32\%& [-] & 5 & [-] \\ \hline
\end{tabular} 
\caption{Optimization results for the Branin function over 10 repetitions}
\label{TableBranin}
\end{table}

\begin{figure}[H]
\hspace{-0.5cm}
 \includegraphics[width=1.1\linewidth]{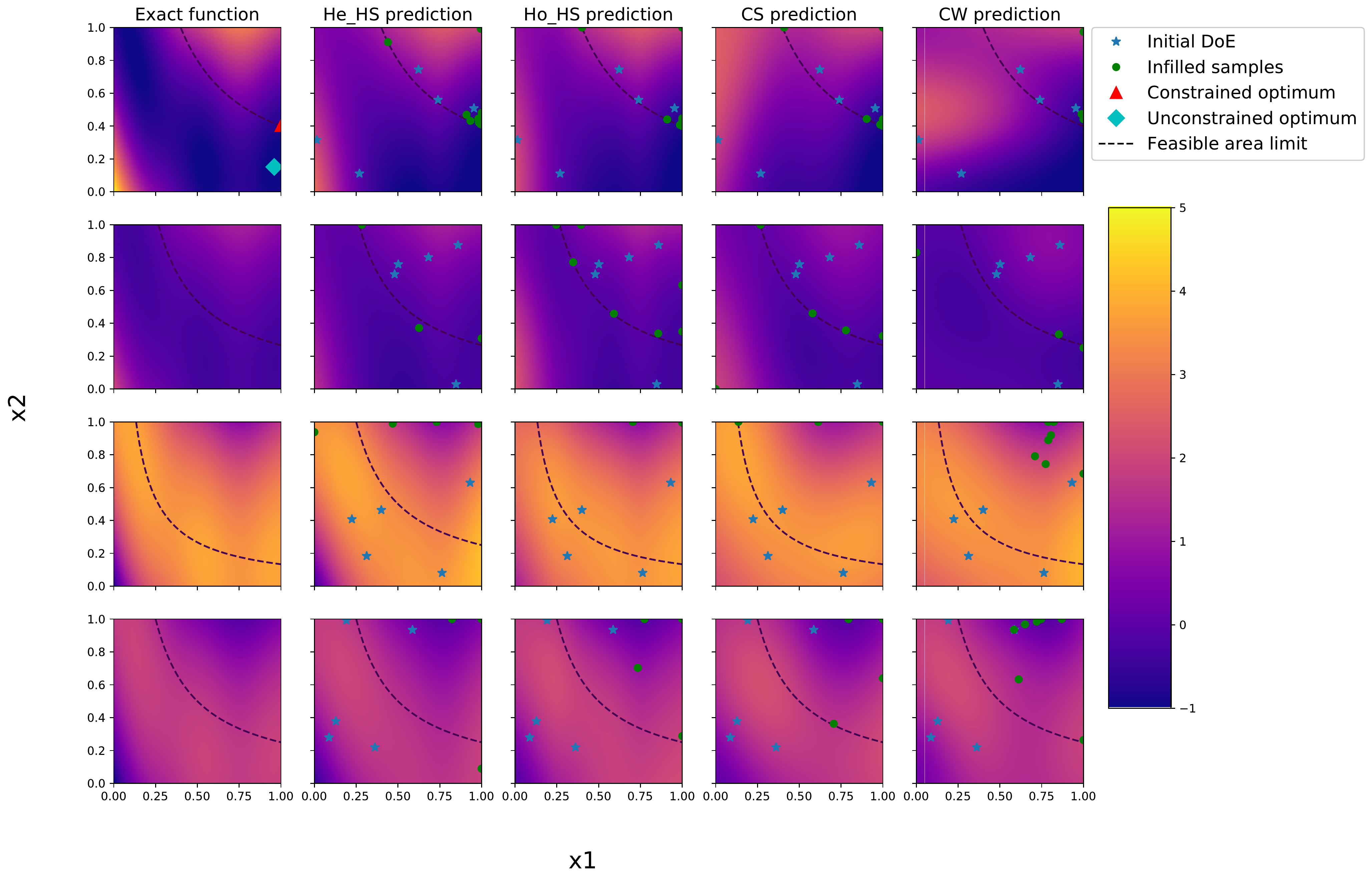}
  \caption{Repartition of the infilled samples as a function of the considered methods (horizontally) and the problem categories (vertically). The optimization techniques rely on heteroscedastic dimension-wise decomposition (He$\_$HS),  homoscedastic dimension-wise decomposition (Ho$\_$HS), compound symmetry decomposition (CS) and category-wise separate surrogate modeling (CW)}
  \label{BraninInfilledPts}
  \end{figure}

For this first analytical benchmark, the proposed mixed variable adaptations of EGO perform better than the standard category-wise EGO and the GA. More specifically, it can be noted that they provide a more robust convergence towards the neighborhood of the problem optimum. On the contrary, both the GA and the category-wise EGO are not provided with enough iterations to consistently identify the area of interest of the search space, most notably for the GA which does not converge to the problem optimum neighborhood in any of the repetitions. Due to the simplicity of this first analytical test-case, no noticeable difference in performance between the proposed mixed variable EGO adaptations can be noticed. For illustrative purposes, the data samples infilled by the presented methods during one of the repetitions as well as the resulting objective function predictions (at the end of the optimization process) are presented in Figure \ref{BraninInfilledPts}.

As can be expected, it can be noticed that the mixed variable SMBDO techniques identify the area of interest of the problem after few iterations so that a large part of the infilled samples are within the neighborhood of the problem optimum. The category-wise EGO, instead, needs to explore the mixed variable search space before being able to identify the correct category as well as the area of interest. By consequence, a larger amount of function evaluations are performed in the non-optimal categories and areas of the problem. Furthermore, it can also be noticed that the category-wise surrogate modeling result in a more inaccurate surrogate modeling of the objective function, especially for the first two categories of the problem, due to the fact that it relies on a smaller amount of data in order to characterize the GP. 
\newpage

\subsection{Augmented Branin function}
\begin{figure}[h!]
 \includegraphics[width=0.9\linewidth]{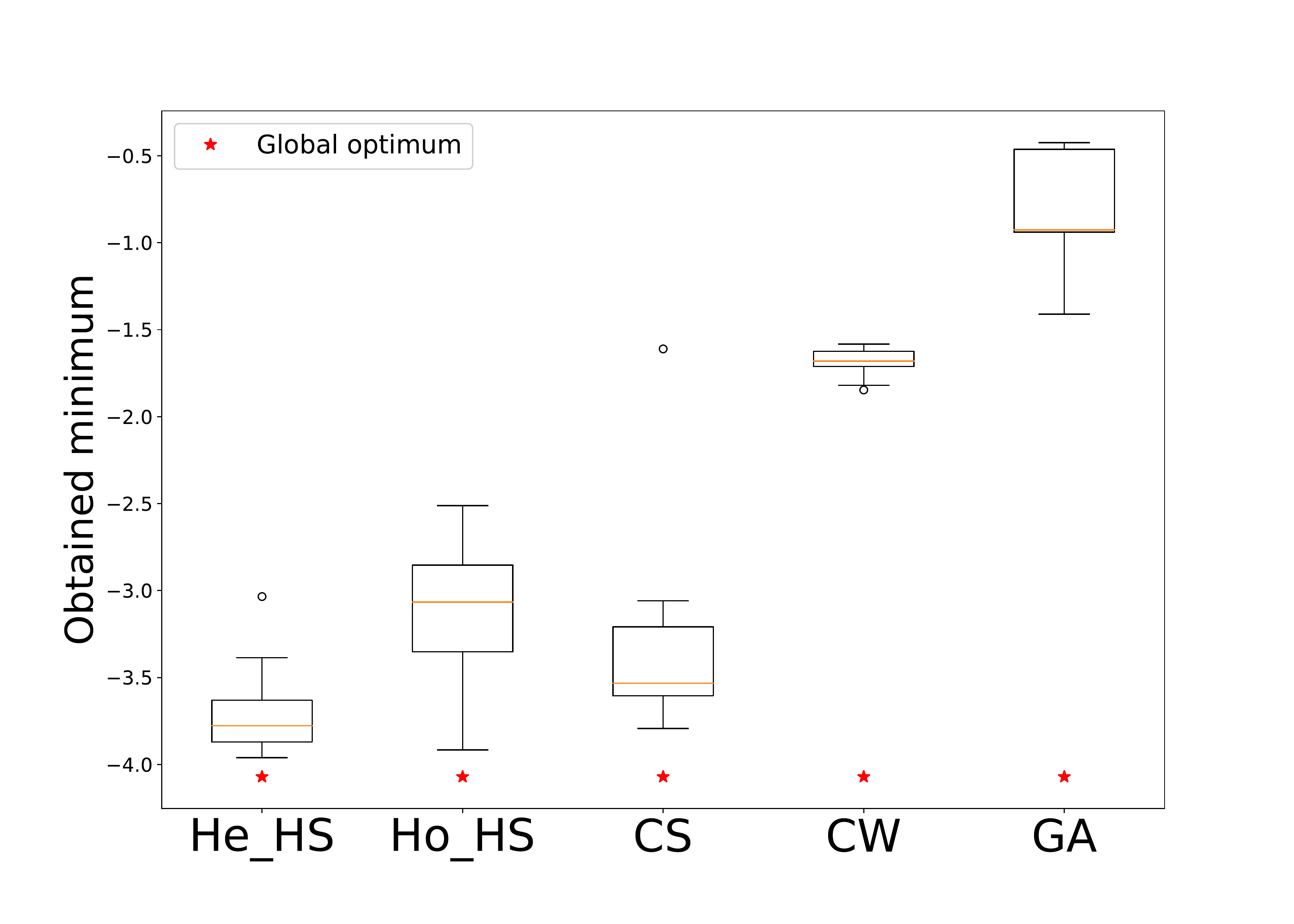}
  \caption{  Optimization results for an augmented Branin function obtained over 10 repetitions}
  \label{BraninAugRes1}
  \end{figure}

\begin{figure}[h!]
 \includegraphics[width=0.9\linewidth]{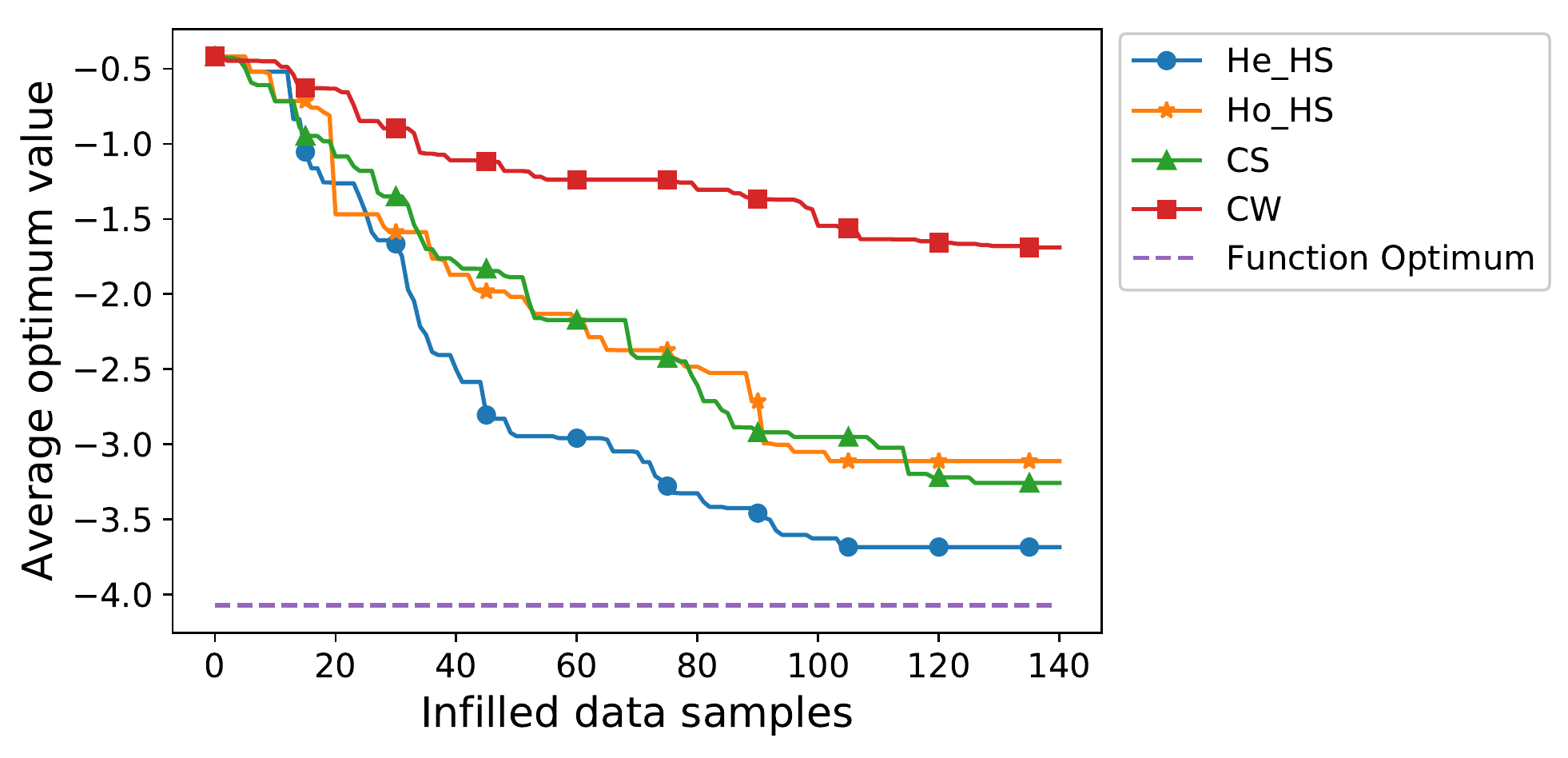}
  \caption{Evolution of the obtained optimum value for the optimization of the augmented Branin test-case as a function of the number of infilled data samples averaged over 10 repetitions}
  \label{BraninAugRes2}
  \end{figure}

In order to assess the influence of the number of continuous dimensions on the performance of the proposed optimization techniques, an augmented version of the previously described Branin function with 10 continuous decision variables and similar discrete categories as the previous test-case  is considered. This second analytical benchmark problem is defined as follows:
\begin{eqnarray}
&& \min   f(x_1,...,x_{10},z_1,z_2)   \\ 
&& \text{w.r.t.} \quad x_1...,x_{10},z_1,z_2 \nonumber
\end{eqnarray}
\begin{equation}
\text{s.t.} \quad g(x_1,...,x_{10},z_1,z_2) \geq 0
\end{equation}
with:
\begin{eqnarray*}
x_1, ..., x_{10} \in [0,1], \quad z_1 \in \{0,1\}, \quad z_2 \in \{0,1\} 
\end{eqnarray*}

\noindent where:
\begin{equation}
f(x_1,...,x_{10},z_1,z_2) = \sum_{i = 1}^{5}  s(x_{2i-1},x_{2i},z_1,z_2)
\end{equation}

\begin{equation}
  s(x_i,x_j,z_1,z_2)=\begin{cases}
               h(x_i,x_j) \hfill   \qquad  \mbox{ if }  z_1 = 0  \mbox{ and } z_2 = 0\\
               0.4 h(x_i,x_j) \hfill  \qquad  \mbox{ if } z_1 = 0 \mbox{ and } z_2 = 1\\
               -0.75h(x_i,x_j) +3.0  \hfill  \qquad  \mbox{ if } z_1 = 1 \mbox{ and } z_2 = 0 \\
               -0.5h(x_i,x_j)  +1.4 \hfill  \qquad  \mbox{ if } z_1 = 1 \mbox{ and } z_2 = 1
            \end{cases}
\end{equation}

\begin{equation}
\hspace{-1.5cm}
\begin{aligned}
h(x_i,x_j) = & \left[ \left( (15x_j-\frac{5}{4\pi^2}(15x_i-5)^2+\frac{5}{\pi}(15x_i-5)-6)^2 + \right. \right.
\\ & \left. \left. 10(1-\frac{1}{8\pi})\cos(15x_i-5)+10 \right) -54.8104 \right] \frac{1}{51.9496}
\end{aligned}
\end{equation}

\begin{equation}
g(x_1,...,x_{10},z_1,z_2) = \sum_{i = 1}^{5}  u(x_{2i-1},x_{2i},z_1,z_2)
\end{equation}
\begin{equation}
  u(x_i,x_j,z_1,z_2)=\begin{cases}
               x_ix_j-0.4 \hfill   \qquad  \mbox{ if }  z_1 = 0  \mbox{ and } z_2 = 0\\
               1.5x_ix_j-0.4 \hfill  \qquad  \mbox{ if } z_1 = 0 \mbox{ and } z_2 = 1\\
              1.5x_ix_j-0.2 \hfill  \qquad  \mbox{ if } z_1 = 1 \mbox{ and } z_2 = 0 \\
               1.2x_ix_j-0.3 \hfill  \qquad  \mbox{ if } z_1 = 1 \mbox{ and } z_2 = 1
            \end{cases}
\end{equation}
For the SMBDO techniques, an initial training data set of 60 samples is used and subsequently 140 additional data points are infilled during the optimization process. The GA is initialized with a population of 10 individuals which evolve for 20 generations. The results obtained for the optimization of the augmented Branin function averaged over 10 repetitions are shown in Figure \ref{BraninAugRes1}  and \ref{BraninAugRes2} as well as in Table \ref{TableBranin10}.
\begin{table}
\footnotesize
\begin{tabular}{|c|c|c|c|c|}
\hline
\multirow{ 2}{*}{Method} & Average & Average   & \# of optima in & \# of \\ 
 & optimum value & constraint value  &  the correct category & hyperparameters \\ \hline
He$\_$HS & -3.683 $\pm$ 7.765 \% & 0.058 & 10 & 26 \\ \hline
Ho$\_$HS & -3.112 $\pm$ 12.242 \% & 0.045  & 10   & 22 \\ \hline
CS &  -3.276 $\pm$ 19.102 \% & 0.051 & 9 & 24\\ \hline
CW & -1.690 $\pm$ 5.147 \%& 0.063 & 0  & 10/category\\ \hline
GA & -0.947 $\pm$ 35.908 \% & [] & 2 & [-]\\ \hline
\end{tabular} 
\caption{Optimization results for the augmented Branin function over 10 repetitions }
\label{TableBranin10}
\end{table}
Similarly to the previous test-case, the proposed mixed variable EGO adaptations provide more robust results when few function evaluations are allowed if compared to the reference methods. Furthermore, a slightly better performance of the heteroscedastic hypersphere decomposition based EGO with respect to the other proposed methods can be noticed. However, differently than with the two-dimensional Branin function, it can also be noted that over the 10 repetitions none of the compared methods converges to the actual problem minimum. It can therefore be stated that a larger number of functions evaluations would be required for any of the considered methods to converge to the problem optimum. 


\newpage 

\subsection{Goldstein function}
In order to assess the dependency of the proposed algorithm performance on the discrete dimension of the optimization problem, the last analytical benchmark function to be considered is a modified version of the Goldstein function characterized by 2 continuous variables and 2 discrete variables, each one presenting 3 levels, thus resulting in a total of 9 discrete categories. This analytical test-case is introduced in order to asses the dependency of the various surrogate models performance on the number of categories that characterize the problem. The Goldstein function optimization problem is defined as follows:

\begin{eqnarray}
&& \min   f(x_1,x_2,z_1,z_2)   \\ 
&& \text{w.r.t.} \quad x_1,x_2,z_1,z_2 \nonumber
\end{eqnarray}
\begin{equation}
\text{s.t.} \quad g(x_1,x_2,z_1,z_2) \geq 0
\end{equation}
with:
\begin{eqnarray*}
x_1 \in [0,100], \quad x_2 \in [0,100], \quad z_1 \in \{0,1,2\}, \quad z_2 \in \{0,1,2\} 
\end{eqnarray*}

\noindent where:
\begin{equation}
f(x_1,x_2,z_1,z_2) = h(x_1,x_2,x_3,x_4) 
\end{equation}

\begin{equation}
\begin{aligned}
    h(x_1,x_2,x_3,x_4)  = \quad & 53.3108+ 0.184901x_1-5.02914 x_1^3 \cdot 10^{-6} + 7.72522x_1^4 \cdot 10^{-8}-  \\
        & 0.0870775x_2 - 0.106959 x_3 +  7.98772x_3^3 \cdot 10^{-6} +  \\
     & 0.00242482 x_4 +  1.32851 x_4^3 \cdot 10^{-6} - 0.00146393 x_1 x_2 -  \\
       & 0.00301588 x_1 x_3 - 0.00272291 x_1 x_4+ 0.0017004 x_2 x_3 +  \\ 
       & 0.0038428 x_2 x_4 -  0.000198969 x_3x_4 +  1.86025 x_1x_2x_3 \cdot 10^{-5} -  \\
       & 1.88719 x_1 x_2 x_4 \cdot 10^{-6}+  2.50923 x_1 x_3 x_4 \cdot 10^{-5} - \\ & 5.62199 x_2 x_3 x_4 \cdot 10^{-5}
        \end{aligned}
\end{equation}

\begin{equation}
  g(x_1,x_2,z_1,z_2)= c_1\sin \left(\frac{x_1}{10} \right)^3 + c_2\cos\left(\frac{x_2}{20}\right)^2
\end{equation}

\begin{figure}[h!]
 \includegraphics[width=0.9\linewidth]{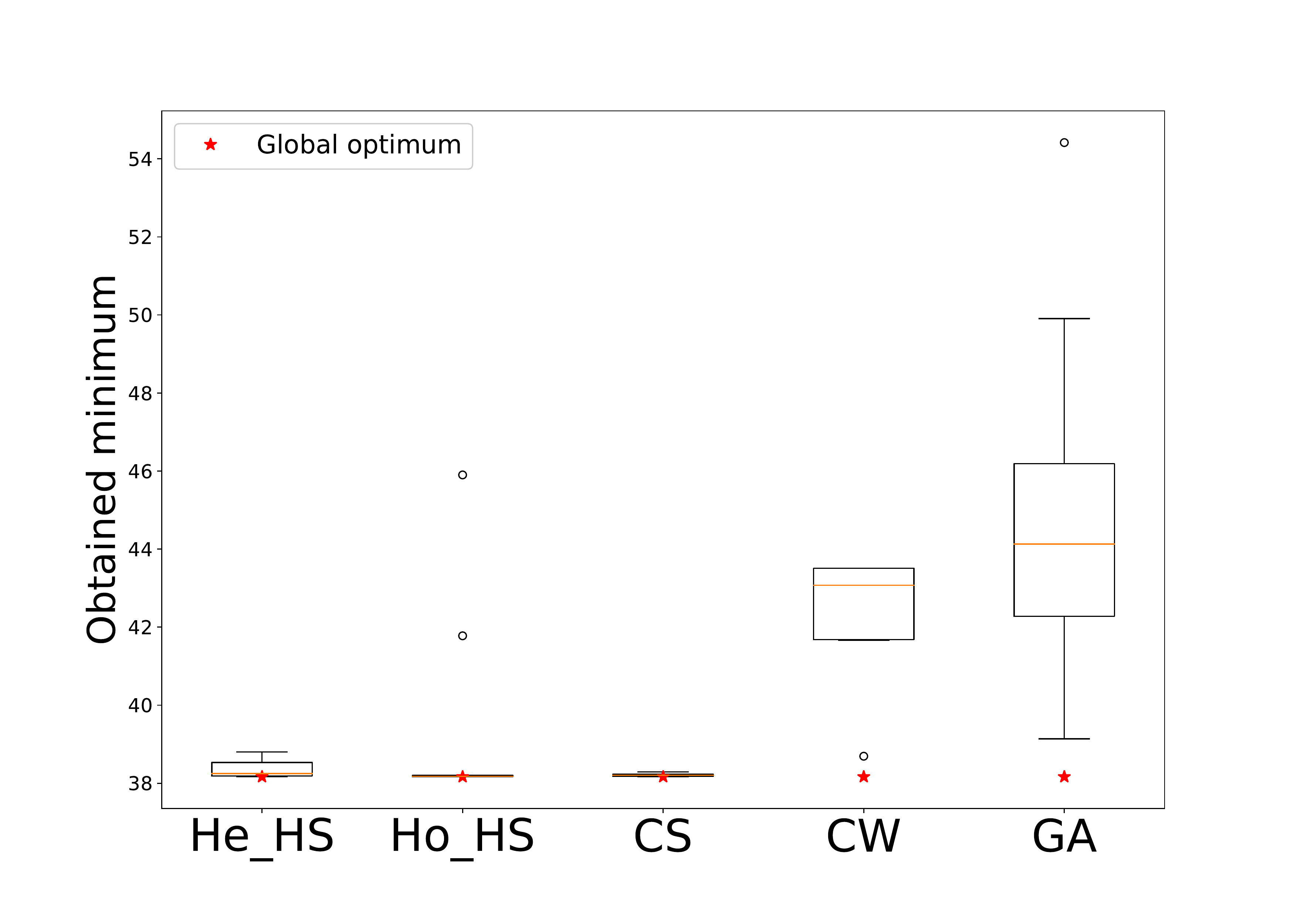}
  \caption{Optimization results for a constrained Goldstein function obtained over 10 repetitions. The used optimization techniques rely on heteroscedastic dimension-wise decomposition (He$\_$HS),  homoscedastic dimension-wise decomposition (Ho$\_$HS), compound symmetry decomposition (CS), category-wise separate surrogate modeling (CW) and a GA}
  \label{GoldsteinRes1}
  \end{figure}

\begin{figure}[h!]
 \includegraphics[width=0.9\linewidth]{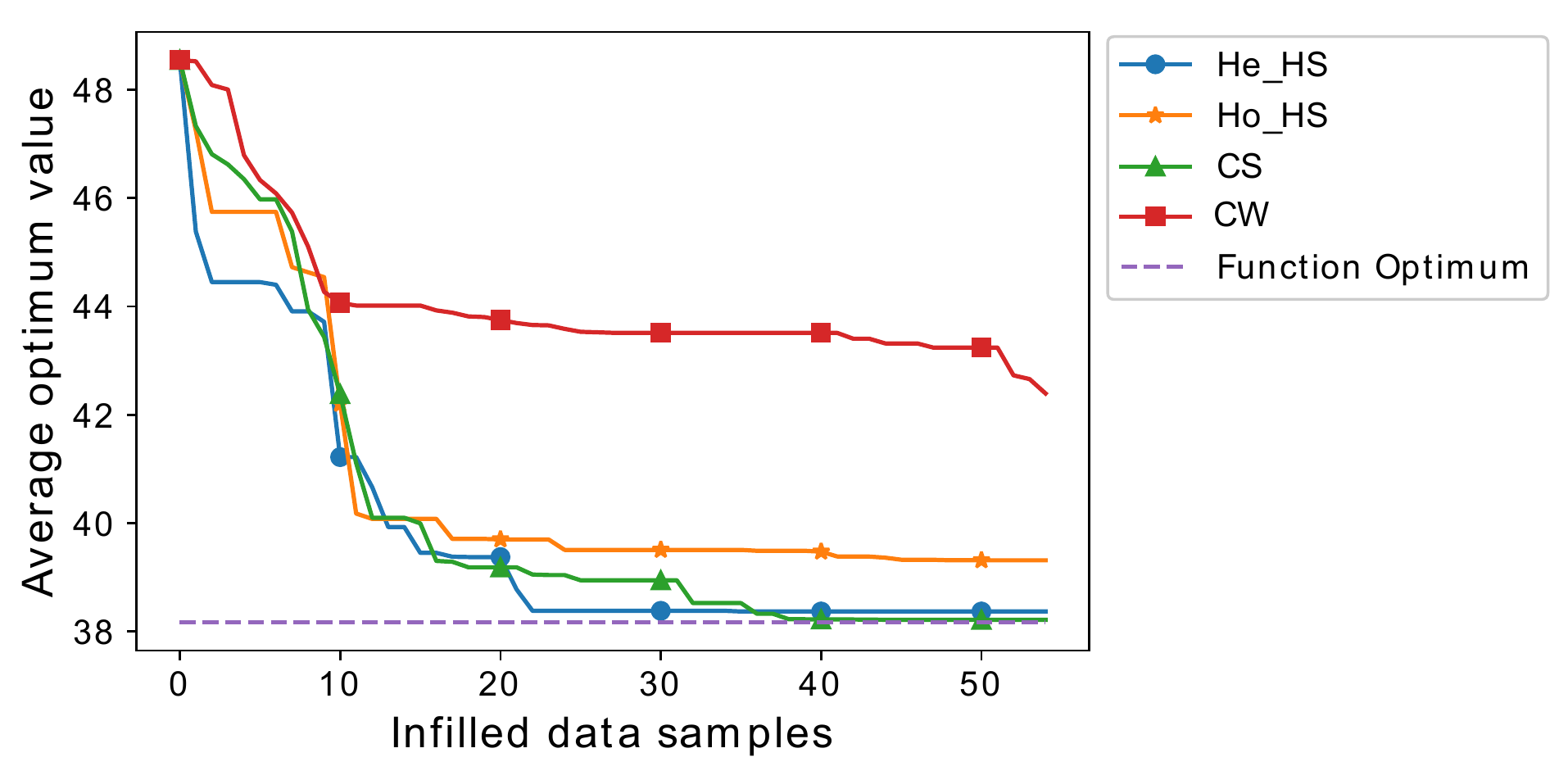}
  \caption{Evolution of the obtained optimum value for the optimization of the Goldstein test-case as a function of the number of infilled data samples averaged over 10 repetitions}
  \label{GoldsteinRes2}
  \end{figure}

\begin{figure}[h!]
\hspace{-0.43cm}
 \includegraphics[width=1.02\linewidth]{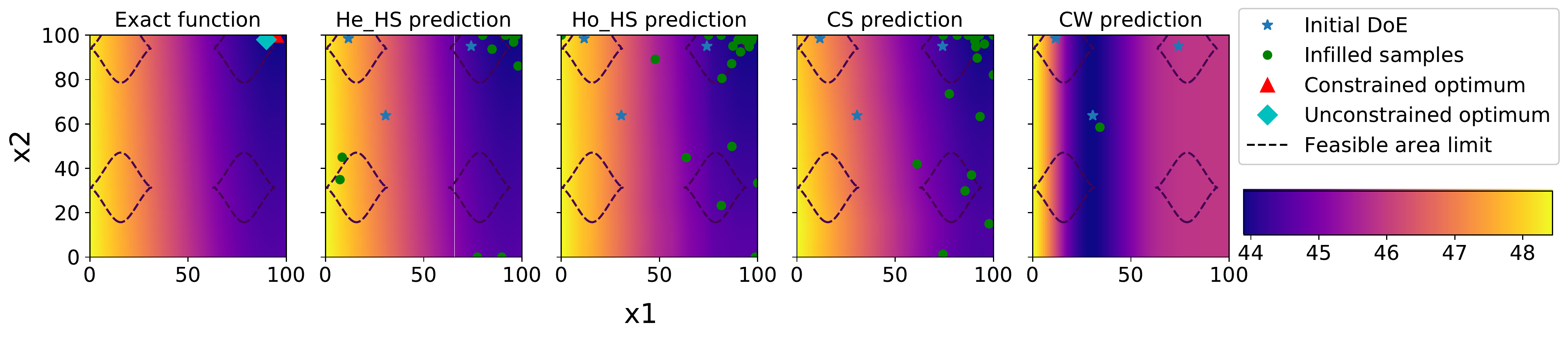}
  \caption{Repartition of the infilled samples as a function of the considered methods. The optimization techniques rely on heteroscedastic dimension-wise decomposition (He$\_$HS),  homoscedastic dimension-wise decomposition (Ho$\_$HS), compound symmetry decomposition (CS) and category-wise separate surrogate modeling (CW)}
  \label{GoldsteinInfilledPts}
  \end{figure}

\noindent The values of $x_3$, $x_4$, $c_1$ and $c_2$ are determined as a function of $z_1$ and $z_2$ according to the relations provided in Table \ref{Gold}.
\begin{table}[h!]
\centering
\begin{tabular}{|c|c c|c c|c c|}
\hline
& \multicolumn{2}{c|}{ $z_1 = 0 $} & \multicolumn{2}{c|}{ $z_1 = 1 $} & \multicolumn{2}{c|}{ $z_1 = 2 $} \\ \hline
\multirow{2}{*}{$ z_2 = 0$} & $x_3 =  20 $&$ x_4 = 20$ & $ x_3 =  50 $&$ x_4 = 20$  & $x_3 =  80$&$ x_4 = 20 $ \\
& $c_1 = 2 $&$ c_2 = 0.5$ & $c_1 = -2 $&$ c_2 = 0.5$  & $c_1 = 1 $&$ c_2 = 0.5$ \\ \hline
\multirow{2}{*}{$ z_2 = 1$} & $x_3 =  20 $&$ x_4 = 50$ & $ x_3 =  50 $&$ x_4 = 50$  & $x_3 =  80 $&$ x_4 = 50 $ \\
& $c_1 = 2 $&$ c_2 = -1$ & $c_1 = -2 $&$ c_2 = -1$  & $c_1 = 1 $&$ c_2 = -1$ \\ \hline
\multirow{2}{*}{$ z_2 = 2$} & $x_3 =  20 $&$ x_4 = 80$ & $ x_3 =  50 $&$ x_4 = 80$  & $x_3 =  80 $&$ x_4 = 80 $ \\
& $c_1 = 2 $&$ c_2 = -2$ & $c_1 = -2 $&$ c_2 = -2$  & $c_1 = 1 $&$ c_2 = -2$ \\ \hline
\end{tabular}
\caption{Characterization of the Goldstein function discrete categories}
\label{Gold}
\end{table}
\noindent For the SMBDO techniques, an initial training data set of 27 samples is used and subsequently 54 additional data points are infilled during the optimization process. The GA is initialized with a population of 8 individuals which evolve for 9 generations. The results obtained for the optimization of the constrained Goldstein function averaged over 10 repetitions are shown in Figures \ref{GoldsteinRes1} and \ref{GoldsteinRes2} as well as in Table \ref{TableGoldstein}. 
\begin{table}
\footnotesize
\begin{tabular}{|c|c|c|c|c|}
\hline
\multirow{ 2}{*}{Method} & Average & Average   & \# of optima in & \# of \\ 
 & optimum value & constraint value  &  the correct category & hyperparameters \\ \hline
He$\_$HS & 38.367 $\pm$ 0.677 \% & 0.349 & 10 & 16\\ \hline
Ho$\_$HS & 39.312 $\pm$ 6.700 \% & 0.136 & 10 & 10\\ \hline
CS & 38.214 $\pm$ 0.130 \%  & 0.096 & 10 & 8 \\ \hline
CW & 42.392 $\pm$  3.401 \% & 0.873 & 1 &  4/category\\ \hline
GA & 44.840 $\pm$ 7.928 \% & [-] & 3 & [-]\\ \hline
\end{tabular} 
\caption{Optimization results for the Goldstein function over 10 repetitions}
\label{TableGoldstein}
\end{table}
In this case, the difference in performance between the category-wise EGO and the proposed mixed variable adaptations is considerable. In fact, it can be seen that said adaptations provide a nearly constant convergence towards the neighborhood of the problem optimum, which is not the case for either the category-wise EGO or the GA. This can be explained by the fact that, when compared with the two variants of the Branin function previously considered, the number of discrete categories characterizing this problem is larger and by consequence the information that a category-wise GP can rely on is very limited. On the other hand, mixed variable GP presented in this article are able to use information provided by data samples in each category of the problem and require therefore less function evaluations in order to converge to the problem optimum when compared to commonly used optimization algorithms. For illustrative purposes, the data points infilled in the optimum category of the problem by the various considered optimization methods during the last repetition are shown in Figure \ref{GoldsteinInfilledPts}. It can be noted that the category-wise EGO is not able to identify the problem category containing the optimum and only one data sample is infilled. On the contrary, the SMBDO techniques based on mixed variable GP are able to quickly identify the optimum category, which results in the majority of the infilled data samples being evaluated in said category.

\subsection{Analytical test-cases result analysis}
By analyzing the results obtained on the previous analytical test-cases, it is shown that by relying on mixed variable GP rather than several continuous ones, it is possible to obtain better optimization results by using the same number of function evaluations. The difference in performance becomes more noticeable when dealing with problems characterized by larger numbers of categories, as can be seen by comparing the results obtained on the two variants of the Branin function (presenting 4 categories) with the ones obtained on the Goldstein function (presenting 9 categories). As mentioned in the previous sections, this difference can be explained by the fact that the larger the number of categories is, the less data is available for the creation of the separate continuous GP used by the category-wise EGO.  Additionally, it is interesting to note that the relative performance of the three mixed variable kernel parameterizations compared in this paper varies depending on the considered test-case. In general, the CS decomposition EGO is more suitable for optimization problems characterized by simple trends and/or large discrete dimensions, whereas the heteroscedastic hypersphere decomposition EGO is more suitable for problems presenting complex trends, but tends to scale poorly with the discrete dimension of the problem. Finally, the  homoscedastic hypersphere decomposition EGO represents a middle-ground between the two previous optimization methods.

\newpage 
\subsection{Launch vehicle propulsion performance}
\begin{figure}[h!]
\centering
 \includegraphics[width=0.9\linewidth]{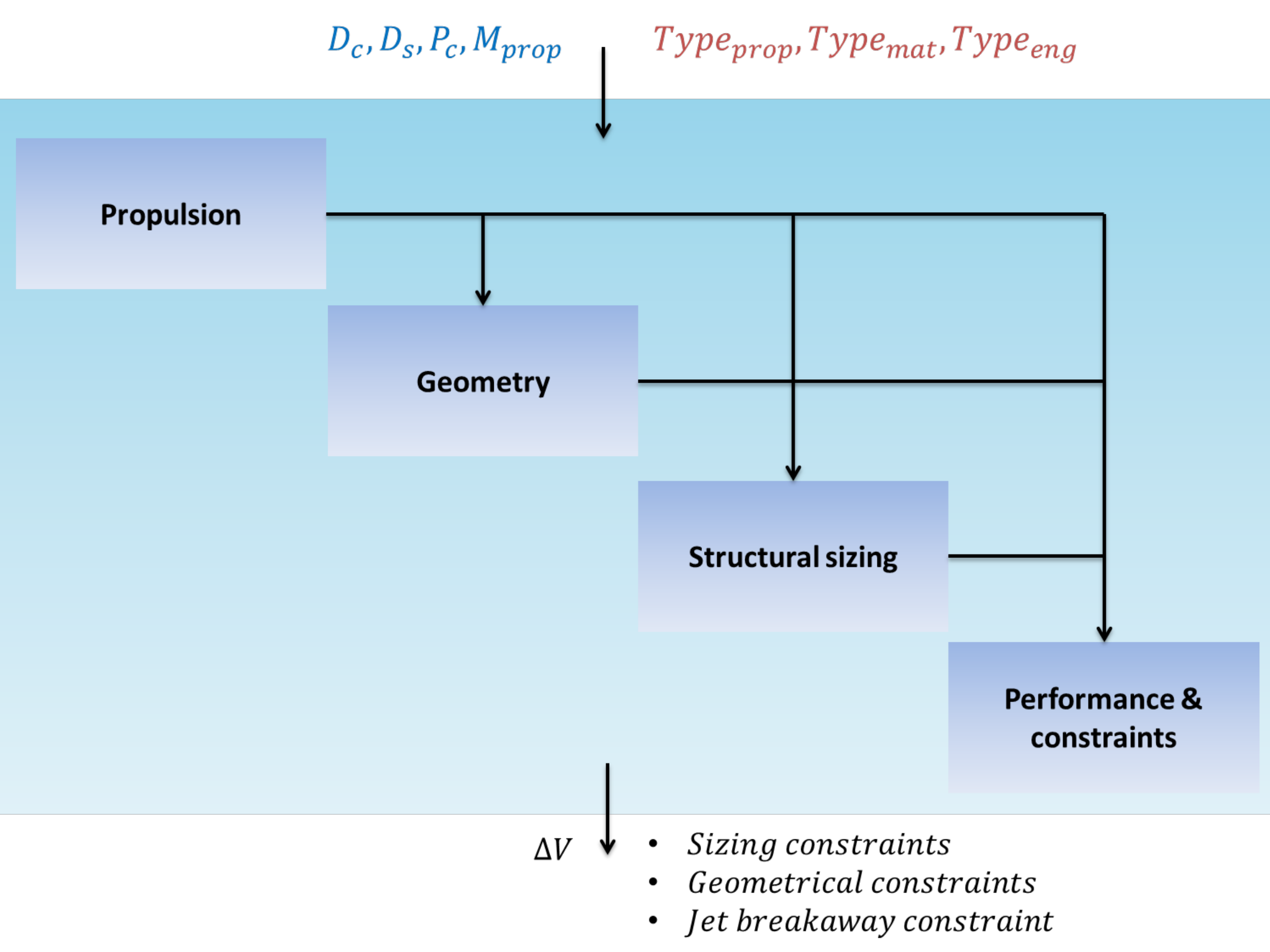}
  \caption{Multidisciplinary design analysis for a solid rocket engine}
  \label{MDO}
  \end{figure}

In order to show the potential applications of the presented mixed variable SMBDO technique for actual engineering design problem, the last test-case considered in this paper is the optimization of a solid rocket engine for a sounding rocket. Sounding rockets carry scientific experiments into space along a parabolic trajectory. Their overall time in space is brief and the cost factor makes them an interesting alternative to heavier launch vehicles as they are sometimes more appropriate to successfully carry out a scientific mission and are less complex to design. The objective of the considered optimization problem is to maximize the speed increment ($\Delta V$) under the geometrical, propulsion and structural constraints.  Three disciplines are involved in the considered test case: the propulsion, the mass budget and geometry design and the structural sizing (Figure \ref{MDO}). 
This mixed variable optimization is characterized by a total of 24 categories. The various continuous and discrete design variables characterizing its performance and constraint functions are detailed in Table \ref{VarBooster}.
\begin{table}
\footnotesize
\begin{tabular}{|c|c|c|c|c|}
\hline
Variable & Nature & Min & Max & Levels \\ \hline
Dc - Nozzle throat diameter [m] & continuous   & 0.2 & 1 & [-] \\ \hline
Ds - Nozzle exit diameter [m] & continuous   & 0.5 & 1.2 & [-]\\ \hline
Pc – Chamber pressure  [bar]  & continuous  & 5 & 300 & [-]\\ \hline
Mprop – Propellant mass  [kg] & continuous  & 2000 & 15000 & [-]\\ \hline
\multirow{ 2}{*}{Type$\_$prop - Type of propellant} & \multirow{ 2}{*}{discrete} &  \multirow{ 2}{*}{[-]} &  \multirow{ 2}{*}{[-]} & Butalite, Butalane, \\ 
& & & &  Nitramite, pAIM-120 \\ \hline
Type$\_$mat - Type of material & discrete  & [-] & [-] & Aluminium, Steel \\ \hline
Type$\_$eng - Type of engine & discrete  & [-] &  [-]&  type 1, type 2, type 3\\ \hline
\end{tabular}
\caption{Variables characterizing the solid rocket engine test-case}
\label{VarBooster}
\end{table}

\begin{figure}[h!]
\centering
 \includegraphics[width=0.5\linewidth]{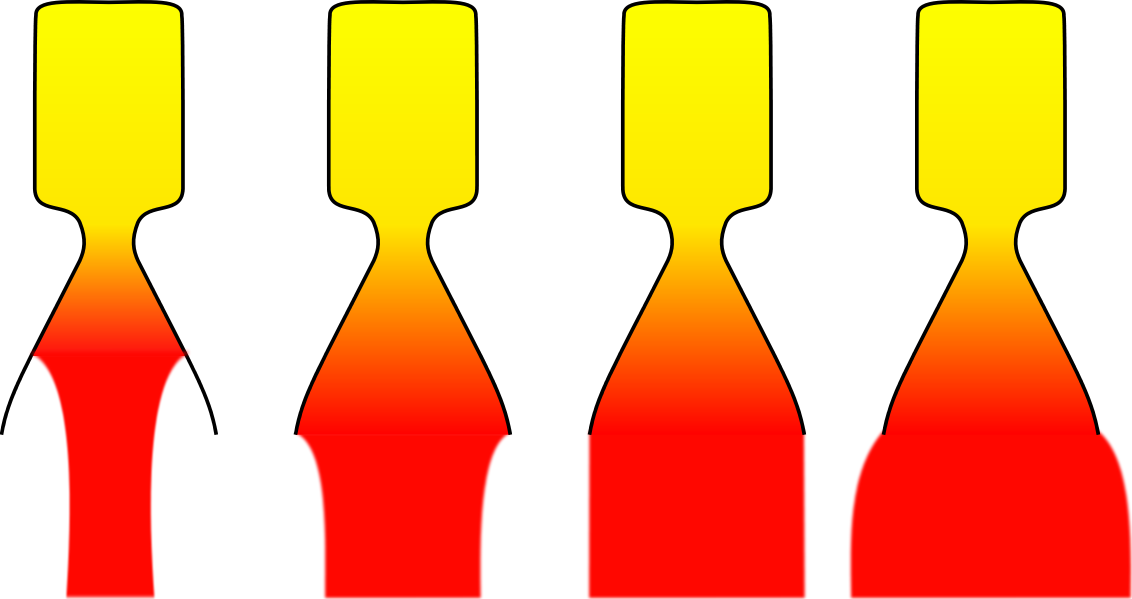}
  \caption{Under or over gas expansion at the nozzle exit}
  \label{nozzle}
  \end{figure}

Different organizations of the solid propellant in the cylindrical motor case may be defined depending on the grain geometry inside the case, which has consequences on the level of thrust of the solid rocket engine along the trajectory. A star grain configuration is used in the simulation as it presents the advantage of providing a constant propellant burning surface along the trajectory and by consequence a constant thrust. Three different engine options, each one with different efficiency and geometrical/physical properties, are considered. In the propulsion discipline, the propellant burning rate is determined as a function of the type of propellant (four options are considered: Butalite, Butalane, Nitramite and pAIM-120) as well as the combustion chamber pressure. Moreover, depending on the mass of propellant and the nozzle geometry (nozzle throat diameter and nozzle exit diameter), the specific impulse and the thrust may be estimated. To avoid jet breakaway (schematically represented in Figure \ref{nozzle}), a constraint on the gas expansion is considered.
In the mass budget and geometry design discipline, the masses of the solid rocket engine case and the nozzle are estimated based on the type of material, the nozzle geometry and the propellant masses. Six constraints relative to the available volume for the propellant, the case geometry and the nozzle geometry are present in the problem. Finally, the structure discipline includes a constraint representing the maximal material stress, as the motor case must be able to withstand the maximum combustion chamber pressure under any possible operating condition. The performance and constraints module estimates the speed increment which combines the rocket engine thrust, the specific impulse, the propellant mass and the booster dry mass. 

The sounding rocket optimization test-case described above can be formally defined as:

\begin{align}
 \min & \qquad  -\Delta V(D_c,D_s,P_c,M_{prop},Type_{prop}, Type_{mat}, Type_{eng}) \\ 
\text{w.r.t.} &  \qquad D_c,D_s,P_c,M_{prop},Type_{prop}, Type_{mat}, Type_{eng} \nonumber \\
\text{s.t.:} & \qquad g_i(D_c,D_s,P_c,M_{prop},Type_{prop}, Type_{mat}, Type_{eng}) \leq 0  \quad \mbox{for i } =1,...,8  \nonumber
\end{align}
with:
\begin{eqnarray*}
D_c \in [0.2,1], \quad D_s \in [0.5,1.2], \quad P_C \in [5,300], \quad M_{prop} \in [2000,15000], \\
Type_{prop} \in \{0,1,2,3\},\quad Type_{mat} \in \{0,1\}, \quad Type_{prop} \in \{0,1,2\}
\end{eqnarray*}

\begin{figure}[h!]
 \includegraphics[width=0.9\linewidth]{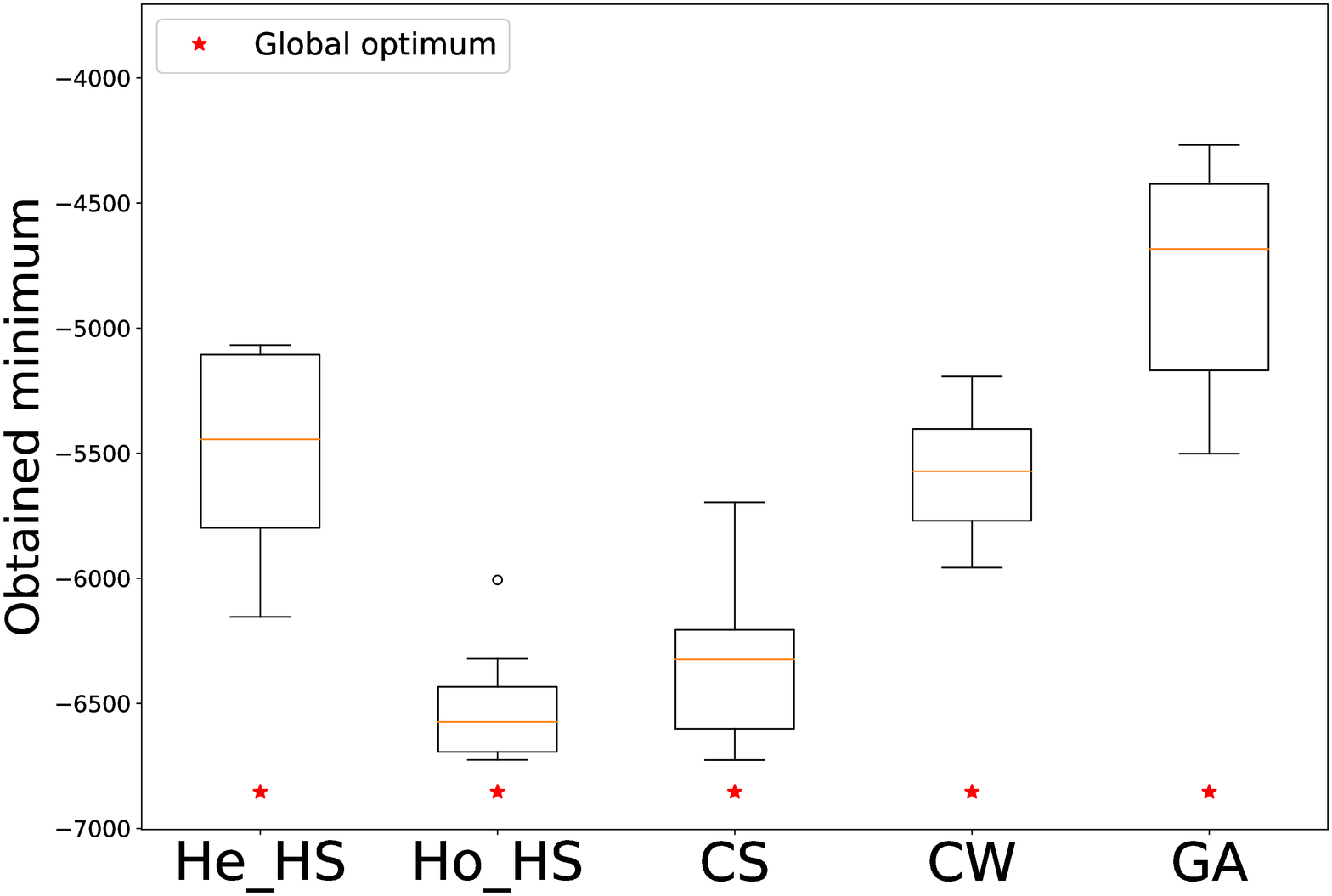}
  \caption{Optimization results for a sounding rocket propulsive performance over 10 repetitions}
  \label{ResBooster1}
  \end{figure}

\begin{figure}[h!]
 \includegraphics[width=0.9\linewidth]{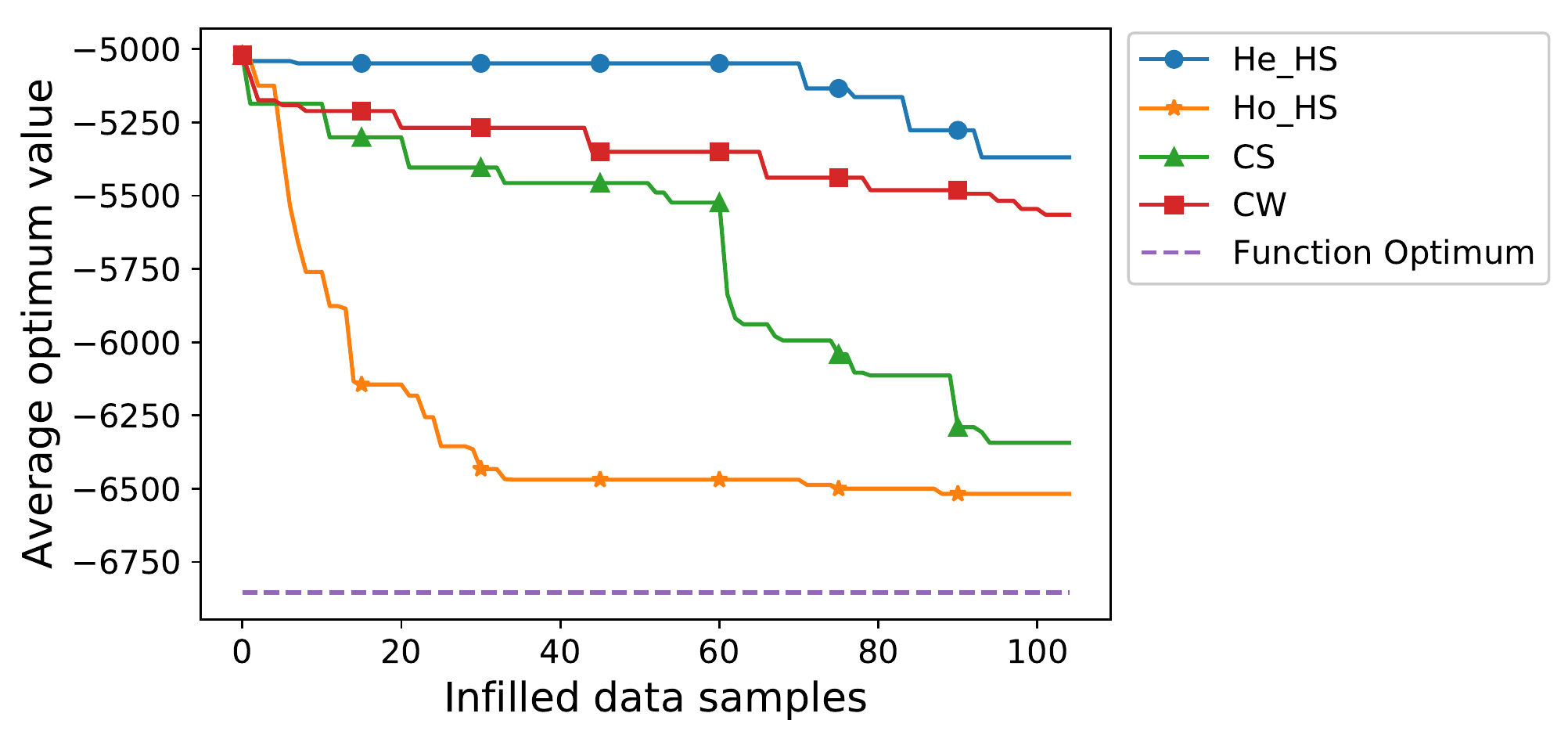}
  \caption{Evolution of the obtained optimum value for the optimization of a sounding rocket propulsive performance as a function of the number of infilled data samples averaged over 10 repetitions}
  \label{ResBooster2}
  \end{figure}

\noindent For the SMBDO techniques, an initial training data set of 96 samples is used and subsequently 104 additional data points are infilled during the optimization process. The GA is initialized with a population of 10 individuals which evolve for 20 generations. The results obtained for the optimization of the sounding rocket engine $\Delta V$ over 10 repetitions are shown in Figures \ref{ResBooster1} and \ref{ResBooster2} as well as in Table \ref{TableBooster}. 
  

\begin{table}
\footnotesize
\begin{tabular}{|c|c|c|c|}
\hline
\multirow{ 2}{*}{Method} & Average    & \# of optima in & \# of \\ 
 & optimum value &  the correct category & hyperparameters \\ \hline
He$\_$HS & -5369.021 $\pm$ 11.545 \%  & 0 & 27 \\ \hline
Ho$\_$HS & -6517.410  $\pm$  3.330  \% & 10 & 18\\ \hline
CS & -6343.094  $\pm$  4.569  \% & 4 & 14\\ \hline
CW & -5564.890  $\pm$ 4.392  \%  & 1 & 8/category\\ \hline
GA & -5177.390  $\pm$ 19.864  \%  & 2 & [-]\\ \hline
\end{tabular} 
\caption{Optimization results for the sounding rocket engine performance over 10 repetitions}
\label{TableBooster}
\end{table}

When compared to the previous analytical test-cases, the optimization of the sounding rocket engine performance is overall more complex, as it is characterized by larger continuous and/or discrete dimensions, a larger number of discrete categories and a larger number of constraints. By analyzing the results obtained for this test-case, it can be noticed that the homoscedastic hypersphere decomposition based EGO performs better than the other mixed-variable SMBDO techniques and than the reference methods, with respect to the average obtained optimum as well as the number of optima in the correct category. The less efficient performance of the CS decomposition EGO can be explained by the fact that the objective and some of the constraint functions trends are complex and therefore difficult to model by relying solely on one covariance and one variance hyperparameter per discrete variable. For what concerns the heteroscedastic hypersphere EGO, due to the large discrete dimension of the problem, the number of hyperparameters required for the characterization of the kernel matrices is too large when compared to the number of data samples that are provided for the training of the surrogate model. Similarly to the augmented Branin test-case presented in the previous section, the number of function evaluations that are provided to the optimization algorithm is not sufficient to consistently converge towards the problem optimum and additional computations would be required in order to be able to use the obtained results within an actual system design framework.

\newpage

\section{Conclusions and perspectives}
In this article, various adaptations of the EGO algorithm for the optimization of constrained problems depending simultaneously on continuous and discrete design variables are proposed, discussed and compared on several test-cases. It is shown that by relying on a single mixed variable surrogate model rather than multiple purely continuous ones, it is possible to reduce the number of functions evaluations necessary to converge towards the optimum value of a given problem. The gain in efficiency is more noticeable for optimization problems  characterized by a large numbers of discrete variable value combinations, as in these cases the number of required separate surrogate models becomes too large with respect to the amount of available data. The proposed optimization algorithms result therefore promising for the optimization of constrained mixed variable optimization problems which present computationally intensive objective and/or constraints functions, as they require fewer evaluations in order to identify the problem optimum. However, it must be noted that the gain in efficiency provided by the presented mixed variable SMBDO techniques tends to diminish when dealing with simpler problems and shorter computation times due to the computational overhead related to the surrogate model training and to the infill criterion optimization.

Three different parameterizations of the discrete part of the GP kernel are considered and compared in this paper. Each parameterization is based on a number of assumptions and is characterized by a number of hyperparameters which scales with the discrete dimension of the problem. Generally speaking, larger numbers of hyperparameters imply more accurate surrogate models. However, it is shown in the presented results that relying on the theoretically most accurate surrogate model does not necessarily yield the best optimization results. In fact, for some of the presented test-cases, the Compound Symmetry decomposition based EGO yields the best optimization results, while being the SMBDO technique relying on a 'simpler' surrogate modeling technique. This can be explained by the fact that large amounts of data are necessary in order to properly identify the optimal hyperparameter values of complex surrogate models, and as the optimizations presented in the previous sections are performed with limited function evaluations, this is not always possible. Simpler surrogate models, on the other hand, provide poor modeling of functions with complex trends and/or correlations, but are usually characterized by fewer hyperparameters which can therefore be correctly optimized with smaller amounts of data. Among the three parameterizations considered in this paper, the one that seems to provide the most robust results and that scales better with both the dimension and the complexity of the problem is the homoscedastic hypersphere decomposition, whereas the compound symmetry decomposition performs less efficiently when dealing with complex trends and the heteroscedastic hypersphere decomposition scales poorly with the discrete dimension of the problem. In general, it can be concluded that in order to choose the most suitable parameterization of the GP mixed variable kernel, a trade-off  between the complexity of the chosen surrogate model and the amount of data that can be provided in order to properly train it is necessary, and that no universally valid parameterization exists.

In order to provide a better assessment of the proposed SMBDO technique performance and scaling capabilities with respect to the continuous and discrete dimensions size as well as problem complexity, further engineering related test-cases from various industry domains should be considered. Finally, it is important to note that the proposed mixed variable optimization algorithm is based on the assumption that all the discrete variables the problem depends on are nominal, which means that no relation of order between their possible values exist. However, this is not always true, as the considered problems can also depend on integer variables, which are inherently ordered, as well as on so-called 'ordinal' discrete variables for which a relation of order can be defined (\emph{e.g.,} small, medium, large type of variables). It could therefore be interesting to explore the possibility of differentiating the discrete kernels between the ones for ordinal and nominal variables, thus allowing to maximize the information that can be extracted from the data samples used for the creation of the surrogate models.


\subsubsection*{Acknowledgements}
This research is co-founded by the Centre National d'\' Etudes Spatiales (CNES) and by the Office National d' \' Etudes et de Recherches Aerospatiales (ONERA - The French Aerospace Lab) within the context of a PhD thesis.

\bibliographystyle{abbrv}
\bibliography{References}
\end{document}